\documentclass[a4paper,10pt]{article}
\usepackage{graphicx}  
\usepackage[colorlinks=true, allcolors=blue]{hyperref}

\usepackage{amsthm}

\usepackage{amsmath,amsfonts,amssymb}
\usepackage{graphicx}
\usepackage{geometry}
\usepackage{etoolbox}
\usepackage[numbers]{natbib}
\AtBeginEnvironment{thebibliography}{\setlength{\itemsep}{0pt}\setlength{\parsep}{0pt}}
\usepackage{hyperref}

\usepackage{xcolor}
\usepackage{enumitem}
\usepackage{amsthm}
\usepackage{comment}
\usepackage{caption}
\usepackage{subcaption}
\usepackage{placeins}

\usepackage{booktabs}
\usepackage{multirow}
\usepackage{tcolorbox}
\usepackage{mathrsfs}
\usepackage{cases}

\numberwithin{equation}{section}
\def\bi{\bm{\mathrm{i}}}

\DeclareMathOperator{\repart}{Re}
\DeclareMathOperator{\impart}{Im}

\newcommand{\wh}[1]{\widehat{#1}}

\newcommand{\inner}[2]{\left\langle {#1}, {#2}\right\rangle}
\newcommand{\norm}[1]{\left\| {#1} \right\|}
\newcommand{\smnorm}[1]{\left| {#1} \right|}

\newcommand{\rmd}{\mathrm{d}}
\newcommand{\bbR}{\mathbb{R}}

\newcommand{\bbC}{\mathbb{C}}

\makeatletter
\newcommand{\Rmnum}[1]{\@Roman{#1}}
\makeatother

\newtheorem{theorem}{Theorem}[section]

\newtheorem{remark}{Remark}
\newtheorem{lemma}{Lemma}[section]

\newtheorem{exm}{Example}[section]

\usepackage{bm}
\usepackage{lineno}
\usepackage{authblk} % new
\newenvironment{keyword}{\noindent\textbf{Keywords: }}{}
\newenvironment{MSC2020}{\noindent\textbf{MSC2020: }}{}

\title{Unconditional optimal convergence of a relaxation Crank-Nicolson difference scheme for 2D space-fractional NLS} 

\author[1]{Lisen Ding\thanks{E-mail: dingmath15@smail.xtu.edu.cn}}
\author[1]{Dongling Wang\thanks{E-mail: wdymath@xtu.edu.cn; Corresponding author.}}
\affil[1]{Hunan Key Laboratory for Computation and Simulation in Science and Engineering, School of Mathematics and Computational Science, Xiangtan University, Xiangtan, Hunan 411105, China.}

\begin{document}

\maketitle

\begin{abstract}
In this paper, we develop a linearly implicit relaxation Crank-Nicolson difference scheme for the two-dimensional
space-fractional nonlinear Schr{\"o}dinger (SFNS) equation, where the fractional Laplacian
is approximated by a simple fractional central difference scheme. The proposed scheme is uniquely solvable and conserves discrete mass for general real-valued nonlinear coefficients, and a modified discrete energy in the cubic case. By introducing a modified error system employing a cut-off technique, we establish the optimal-order error estimates without mesh-ratio restriction.
Specifically, the primary variable achieves second-order temporal and spatial convergence in the discrete
$H_h^{\alpha/2}$-norm, and the relaxation variable satisfies the same convergence in the discrete $L^2$-norm.
Some numerical experiments confirm the relevant theoretical results.
\end{abstract}

\begin{MSC2020}
    35R11, 65M06, 65M12, 35Q55
\end{MSC2020}

\section{Introduction}

Consider the following two-dimensional space-fractional nonlinear Schr{\"o}dinger (SFNS) equation with initial condition:
\begin{alignat}{2}
    \bi u_t - (-\Delta)^{\frac{\alpha}{2}}u + g \left(|u|^2\right)u &= 0,  &\qquad& (x,y) \in \bbR^2, \quad 0<t<T, \label{eq-SFNS} \\
    u(x,y,0) &= u_0(x,y), &\qquad& (x,y) \in \bbR^2, \label{eq-inicd}
\end{alignat}
where  $\bi = \sqrt{-1}$ denotes the imaginary unit, and $g:[0,\infty)\to\mathbb{R}$ is a real-valued continuous function. For example, $g(s) = \lambda_1 s$ or $\lambda_2 s^2$ denotes the cubic or quintic SFNS equation, respectively, where $\lambda_1$ and $\lambda_2$ are constants. 
The fractional Laplacian $(-\Delta)^{\frac{\alpha}{2}}$ is defined by the following hypersingular integral,
\begin{equation}\label{eq-intfraLa}
    (-\Delta)^{\frac{\alpha}{2}}u(\bm{x}) = c_{\alpha} \, \operatorname{p.v.} \int_{\bbR^2} \frac{u(\bm{x})-u(\bm{\zeta})}{|\bm{x}-\bm{\zeta}|^{2+\alpha}} \rmd \bm{\zeta} ,\quad c_{\alpha} = \frac{\alpha 2^{\alpha-1}\Gamma \bigl(1+\frac{\alpha}{2} \bigr)}{\pi \Gamma \bigl(1-\frac{\alpha}{2}\bigr)},
\end{equation}
where $1<\alpha\leq2$ in \eqref{eq-intfraLa}, $\bm{x}=(x,y)$, $\bm{\zeta}=(\zeta_1,\zeta_2)$, ``$\operatorname{p.v.}$'' denotes the Cauchy principal value, and $c_\alpha$ is the normalization constant.

When $\alpha=2$, \eqref{eq-SFNS} reduces to the classical Schr{\"o}dinger equation, which is a basic model in quantum mechanics and can be derived from Feynman's path integrals associated with Brownian trajectories \cite{feynman1965quantum}. Motivated by this connection, Laskin extended the path integral from Brownian trajectories to L{\'e}vy flight paths and developed fractional quantum mechanics \cite{laskin2000fractional}. More specifically, the $\alpha$-stable L{\'e}vy process leads to the space-fractional Schr{\"o}dinger equation. The space-fractional Schr{\"o}dinger equation has found applications in physics, including optical beam propagation \cite{longhi2015fractional,zhang2015propagation}, condensed-matter systems \cite{stickler2013potential}, superfluid hydrodynamics \cite{Tayurskii2012}, and other fields \cite{guo2006some,kirkpatrick2013continuum}. From a mathematical viewpoint, various properties of the SFNS
equation have been extensively studied, such as the well-posedness and global existence of solutions \cite{guo2010global,guo2013well,hong2015fractional} and the existence of ground-state solutions \cite{secchi2013ground,li2019ground}.

Generally speaking, the exact solutions to the space-fractional Schr{\"o}dinger equation are difficult to obtain due to the nonlocality of the fractional Laplacian. Therefore, to investigate the dynamical behavior of solutions, one can use numerical approximations as a practical alternative. In recent years, several numerical methods have been developed for solving the SFNS equation, including finite element methods \cite{li2017galerkin,chen2025two}, finite difference methods \cite{wang2013crank,wang2014linearly,zhao2014fourth,wang2015energy,zhang2019conservative} and spectral methods \cite{zhai2019error,shen2024efficient,jia2026optimal}.
Among these methods, the fractional Laplacian employed in  most studies  is defined as the sum of the Riesz fractional derivatives along the coordinate directions, whose two-dimensional Fourier symbol is $|\theta_1|^\alpha+|\theta_2|^\alpha$. However, the integral fractional Laplacian defined in \eqref{eq-intfraLa} couples different spatial directions, and its corresponding Fourier symbol is $(\theta_1^2+\theta_2^2)^{\alpha/2}$. Although the two operators coincide at $\alpha=2$, they are generally not identical for $1<\alpha<2$. In particular, the non-separability of variables in the integral fractional Laplacian introduces additional difficulties for the numerical discretization of the corresponding space-fractional equations.
For the two-dimensional SFNS equation with the integral fractional Laplacian, Wang et al.\ \cite{wang2022linear} analyzed a conservative three-level linearized difference scheme and obtained second-order convergence in the discrete maximum norm. Hu et al.\ \cite{hu2021convergence} studied a structure-preserving modified Crank-Nicolson difference scheme, proved second-order convergence in the discrete $L^2$-norm using a cut-off technique, and developed a fast implementation. 
While rigorous numerical analysis for fractional equations has been established, the fully discrete analysis of a relaxation scheme with staggered auxiliary variables remains relatively limited, which motivates the present study.

Besse \cite{besse2004relaxation} proposed a linearly implicit relaxation scheme for the classical cubic nonlinear Schr{\"o}dinger equation. The scheme employs staggered auxiliary variables to approximate the nonlinear coefficient and preserves discrete mass and energy. Besse et al. \cite{besse2021energy} subsequently constructed an energy-preserving relaxation scheme for general power-law nonlinearities.
Currently, the relaxation-type method has been successfully extended to SFNS equations. 
Li et al. \cite{li2020relaxation} proposed a relaxation Galerkin FEM for the one-dimensional cubic SFNS equation. Cheng et al. \cite{cheng2020energy} developed an energy-preserving relaxation Fourier spectral method for the one-dimensional SFNS equation with general power-law nonlinearities.
More recently, for the classical equation with a general smooth nonlinearity, Zouraris \cite{zouraris2023error} analyzed a fully discrete relaxation central difference scheme incorporating a half-step initialization and obtained optimal second-order estimates without a mesh-coupling condition. Zhou et al. \cite{zhou2025decoupled} and Gao et al. \cite{gao2025unconditional} established the unconditional optimal-order error estimates for relaxation compact difference schemes for the scalar and coupled nonlinear Schr{\"o}dinger equations, respectively. 
Motivated by these works, we propose a relaxation Crank-Nicolson difference (RCND) scheme for \eqref{eq-SFNS} and analyze its unconditional convergence, where the integral fractional Laplacian \eqref{eq-intfraLa} is approximated by using the second-order fractional central difference operator from \cite{hao2021fractional}. 

In the convergence analysis of the one-dimensional SFNS equation, the fractional Sobolev inequality 
\begin{equation*}
    \|v\|_{\infty} \leq C \|v\|_{H_h^{\frac{\alpha}{2}}}, \qquad \forall v \in H_h^{\frac{\alpha}{2}}
\end{equation*}
is commonly used to prove the boundedness of the numerical solution in the $L^{\infty}$-norm.
However, the above inequality does not hold in two-dimensional space, which presents a significant challenge in the convergence analysis of the two-dimensional SFNS equation. To overcome this difficulty, we establish the following uniform discrete $L^{\infty}$-norm estimate in Lemma \ref{lem-uniform_maxest}:
\begin{equation*}
    \norm{v}_{\infty} \leq C_{\alpha} \left( \norm{v} + \norm{(-\Delta_h)^{\frac{\alpha}{2}} v} \right).
\end{equation*}
With this preparation, we can further analyze the convergence of the proposed RCND scheme. To avoid assuming that the discrete $L^\infty$-norm of the numerical solution is bounded, we introduce a modified error system with two cut-off coefficients, and then establish the closed $L^2$-norm estimates for the auxiliary errors. Then, we prove that both cut-off coefficients are equal to one, so the modified error system is consistent with the original error system. Finally, we can derive the unconditional optimal-order error estimates for the proposed scheme. 

The main contributions of this paper are summarized as follows:
\begin{itemize}
    \item The unique solvability of the RCND scheme is strictly proved. In addition, the RCND scheme preserves  discrete mass conservation for general real-valued function $g$, while also preserving  discrete energy conservation in the  case of cubic nonlinearity. 
    \item In two-dimensional space, we develop a uniform discrete $L^{\infty}$-norm estimate, see Lemma \ref{lem-uniform_maxest}. This estimate plays a crucial role in the subsequent analysis. This estimate plays a crucial role in the subsequent analysis. Finally, we establish the unconditional optimal second-order convergence for the RCND scheme. 
\end{itemize}

The remainder of this paper is organized as follows. In Section \ref{sec-preliminaries}, we  introduce a discrete scheme for the fractional Laplacian and some auxiliary lemmas. In Section \ref{sec-scheme}, we present the relaxation Crank-Nicolson difference scheme and establish its solvability and conservation properties. Next, the unconditional convergence analysis of the proposed scheme is given in Section \ref{sec-convergence}. Section \ref{sec-three-dimensions} extends the method and analysis to the three-dimensional space. Finally, several numerical experiments are presented in Section \ref{sec-numerics}, and a brief conclusion is drawn in Section \ref{sec-conclusion}.

\section{Preliminaries}\label{sec-preliminaries}

In this section, we introduce the basic notations, discrete fractional Sobolev norms, and several auxiliary lemmas used to construct and analyze the numerical method.

For the system \eqref{eq-SFNS}--\eqref{eq-inicd}, it is not feasible to perform numerical simulation in the full space $\bbR^2$. Under the suitable conditions, the authors in \cite{felmer2012positive} proves that when $|\bm{x}| \rightarrow +\infty$, the solution of the SFNS equation decays rapidly to zero.   
Therefore, we restrict the problem to a bounded domain $\Omega = (a,b)^2$ and impose the homogeneous extended Dirichlet boundary condition
\begin{equation}\label{eq-bondcd}
    u(x,y,t) = 0,  \quad  (x,y)\in\Omega^c:=\bbR^2 \backslash \Omega, \quad t\in [0,T],
\end{equation}
where $a \ll 0$ and $b \gg 0$. Therefore, the space $\bbR^2$ in \eqref{eq-SFNS}--\eqref{eq-inicd} is replaced with the bounded domain $\Omega$.

\subsection{Notations}

Let $M\in\mathbb{Z}_+$, set $h = (b-a)/M$, and define the spatial grid points $(x_i,y_j) := (a + ih, a + jh)$ for $0 \leq i,j \leq M$. Then, we define the following grid function spaces: 
\begin{equation*}
    \begin{aligned}
        &\mathsf{C}_h = \left\{ v=\{v_{i, j}\} \,\middle\vert\, v_{i,j}\in \bbC, 0 \leq i,j \leq M\right\}, \quad \mathsf{C}_h^0 = \left\{ v \in \mathsf{C}_h \,\middle\vert\,  v_{0, j}=v_{M, j} = v_{i, 0}=v_{i,M}=0  \right\}, \\
        &\mathsf{R}_h = \left\{ v=\{v_{i, j}\} \,\middle\vert\, v_{i,j}\in \bbR, 0 \leq i,j \leq M\right\}, \quad \mathsf{R}_h^0 = \mathsf{R}_h \cap  \mathsf{C}_h^0 .
    \end{aligned}
\end{equation*}

For  $v,w\in\mathsf C_h$, we define the discrete inner product and the associated norm by
\begin{equation}
    \inner{v}{w} := h^2 \sum_{i=1}^{M-1} \sum_{j=1}^{M-1} v_{i,j} \bar{w}_{i, j}, \quad \norm{v} := \sqrt{\inner{v}{v}}, \quad \norm{v}_\infty := \max_{0 \leq i,j \leq M} |v_{i,j}|,
\end{equation}
where $\bar{w}$ is the complex conjugate of $w$.

For the temporal grid, let $K\in\mathbb{Z}_+$, $\tau = T/K$, and $t_k = k\tau$ for $0 \leq k \leq K$. 
For a temporal grid function $\{v^k\}_{k=0}^{K}$, we define the midpoint average and backward difference by
\begin{equation}
    \wh{v}^{\,k+\frac{1}{2}} := \frac{1}{2}\left(v^{k+1} + v^{k}\right), \quad \delta_t v^{k+1} := \frac{1}{\tau}\left(v^{k+1} - v^{k}\right).
\end{equation}
Moreover, we set $t_{1/4}:=\tau/4$ and $t_{k+1/2}:=t_k + \tau/2$ for $0\leq k\leq K-1$.

For $v\in\mathsf C_h^0$, identified with its zero extension, its semi-discrete Fourier transform can be defined as
\begin{equation*}
    v^{\ast}\left(\theta_1, \theta_2\right) = h^2 \sum_{i=1}^{M-1} \sum_{j=1}^{M-1} v_{i,j} e^{-\bi\left(\theta_1 ih + \theta_2 jh\right)}, \quad \left(\theta_1, \theta_2\right) \in \left[-\frac{\pi}{h}, \frac{\pi}{h}\right]^2,
\end{equation*}
and the corresponding inverse transform is given by
\begin{equation*}
    v_{i,j} = \frac{1}{4\pi^2}\int_{-\frac{\pi}{h}}^{\frac{\pi}{h}}\int_{-\frac{\pi}{h}}^{\frac{\pi}{h}}v^{\ast}\left(\theta_1, \theta_2\right)e^{\bi\left(\theta_1 ih + \theta_2 jh\right)} \rmd \theta_1\rmd\theta_2. 
\end{equation*}
Using Parseval's identity, the frequency-domain representation of the inner product can be expressed as 
\begin{equation*}
    \inner{v}{w} = \frac{1}{4\pi^2} \int_{-\frac{\pi}{h}}^{\frac{\pi}{h}} \int_{-\frac{\pi}{h}}^{\frac{\pi}{h}} v^{\ast}\left(\theta_1, \theta_2\right) \overline{w^{*}\left(\theta_1, \theta_2\right)} \rmd\theta_1 \rmd\theta_2.
\end{equation*}
For $s>0$, define the fractional Sobolev semi-norm as
\begin{equation*}
    |v|_{H_h^s}^2 = \frac{1}{4\pi^2} \int_{-\frac{\pi}{h}}^{\frac{\pi}{h}} \int_{-\frac{\pi}{h}}^{\frac{\pi}{h}} \left(\theta_1^2 + \theta_2^2\right)^s \left|v^{\ast}(\theta_1, \theta_2)\right|^2 \rmd\theta_1 \rmd\theta_2,
\end{equation*}
and the fractional Sobolev norm $\norm{v}_{H_h^s}^2 := \norm{v}^2 + |v|_{H_h^s}^2$. Finally, the fractional Sobolev space is denoted by $H_h^s = \{v \mid \norm{v}_{H_h^s} < \infty\}$.

In what follows, we denote by $C^{\kappa}(\mathbb{R};\mathbb{R})$ the set of all $\kappa$-times continuously differentiable real-valued functions, where $\kappa\in \mathbb{Z}_{+} \cup \{0\}$. 
Two useful lemmas are presented below, which play an important role in the subsequent analysis.

\begin{lemma}[\cite{zouraris2023error}]\label{lem-gfunLip}
Let $\sigma>0$, $\mathsf{R}_{\sigma,h}:=\left\{v\in \mathsf{R}_h:\ \|v\|_{\infty}\leq \sigma\right\}$, and $g\in C^3(\mathbb{R};\mathbb{R})$. Define $g_\sigma^{\kappa}:=\sup_{s\in[-\sigma,\sigma]}|g^{(\kappa)}(s)|$ for $\kappa=1,2$. Then, for any $v^a,\, v^b,\, z^a,\, z^b \in \mathsf R_{\sigma,h}$, the following estimates hold:
\begin{align}
    &\left\|g(v^a)-g(z^a)\right\|\leq g_\sigma^{1} \|v^a-z^a\|, \\
    &\left\|\mathrm{G}\left(v^a,v^b,z^a,z^b\right)\right\| \leq g_\sigma^{2} \|z^a-z^b\|_{\infty}\|v^b-z^b\|     \notag \\
    &\hphantom{{} \left\|\mathrm{G}\left(v^a,v^b,z^a,z^b\right)\right\| \leq} 
    +\left(g_\sigma^{1}+g_\sigma^{2}\|z^a-z^b\|_{\infty}\right)\|v^a-v^b-z^a+z^b\|,
\end{align}
where $\mathrm{G}\left(v^a,v^b,z^a,z^b\right):=g(v^a)-g(v^b)-g(z^a)+g(z^b)$.

\end{lemma}

\begin{lemma}[\cite{zouraris2023error}]\label{lem-realvarLip}
    For $v^a,v^b,z^a,z^b\in\mathsf C_h$, let $\mathrm{S}(v^a,v^b,z^a,z^b) :=|v^a|^2-|v^b|^2-|z^a|^2+|z^b|^2$. Then the following estimates hold:
    \begin{equation}
        \left\||v^a|^2-|z^a|^2\right\| \leq \left(\|v^a\|_\infty+\|z^a\|_\infty\right) \|v^a-z^a\|,
    \end{equation}
    and 
    \begin{equation}
        \begin{aligned} 
            \left\|\mathrm{S}(v^a,v^b,z^a,z^b) \right\| \leq & 2\|z^a-z^b\|_\infty\|v^b-z^b\|\\
            &+\left( \|v^a\|_\infty+\|v^b\|_\infty+\|z^a-z^b\|_\infty \right) \|v^a-v^b-z^a+z^b\|.
        \end{aligned}
    \end{equation}
\end{lemma}

\subsection{Discretization of fractional Laplacian}\label{sec-discrete_fraLa}

A key step in solving the SFNS equation numerically is to discretize the fractional Laplacian \eqref{eq-intfraLa}. 
However, the nonlocality and hypersingularity of \eqref{eq-intfraLa} pose challenges for its discretization. To overcome these challenges, a common strategy is to decompose the integral \eqref{eq-intfraLa} into singular and non-singular parts, and then focus on the elimination of singularity, as detailed in \cite{huang2014numerical,duo2019accurate,minden2020simple,zhang2026high}.
Since the implementation process for these methods is typically complicated, we instead use the second-order fractional central difference method proposed in \cite{hao2021fractional}.
The method is straightforward to implement and its discrete matrix can be assembled efficiently by the fast Fourier transform (FFT).

Let $\mathscr F$ denote the Fourier transform and write $\theta=(\theta_1,\theta_2)$. Define
\begin{equation*}
    \mathcal{W}^{2+\alpha}(\mathbb{R}^2) = \left\{ v \in L^1\left(\bbR^2\right) \,\middle\vert\, \int_{\bbR^2} (1+|\theta|)^{2+\alpha} |\mathscr{F}[v](\theta)| \,\rmd \theta < \infty \right\}.
\end{equation*}
Then, we have the following second-order approximation of the fractional Laplacian.

\begin{lemma}[\cite{hao2021fractional}] \label{lem:discrete_laplacian}
If $v(x,y) \in \mathcal{W}^{2+\alpha}(\mathbb{R}^2)$, then the two-dimensional fractional Laplacian can be approximated by
\begin{equation}
    (-\Delta)^{\frac{\alpha}{2}} v(x,y) = (-\Delta_h)^{\frac{\alpha}{2}} v(x,y) + \mathcal{O}(h^2), \quad 0 < \alpha \leq 2,
\end{equation}
where the discrete fractional Laplacian is defined as
\begin{equation}\label{eq-fulldiscrete-fraLa}
(-\Delta_h)^{\frac{\alpha}{2}} v(x,y) := \frac{1}{h^\alpha} \sum_{i,j \in \mathbb{Z}} a_{i,j}^{(\alpha)} v\left(x+ih, y+jh \right),
\end{equation}
and the Fourier coefficients are given by 
\begin{equation}\label{eq-compute-intcoeff}
a_{i,j}^{(\alpha)} = \frac{1}{4\pi^2} \int_{-\pi}^{\pi} \int_{-\pi}^{\pi} \varphi(\theta)^{\frac{\alpha}{2}} e^{-\bi(i\theta_1 + j\theta_2)} \rmd\theta_1 \rmd\theta_2, \quad \varphi(\theta) = 4 \sin^2\left(\frac{\theta_1}{2}\right) + 4 \sin^2\left(\frac{\theta_2}{2}\right).
\end{equation}
\end{lemma}

In what follows, for any $v\in \mathsf{C}_h^0$, we rearrange it into a column vector, i.e.,
\begin{equation}\label{eq-vcolumn}
    v = \left(v_{1},v_{2},\cdots,v_{M-1}\right)^\top, \quad v_j=\left(v_{1,j},v_{2,j},\cdots,v_{M-1,j}\right), \quad j=1,2,\cdots,M-1.
\end{equation} 
After imposing the boundary condition \eqref{eq-bondcd}, the discrete operator can be written in a matrix-vector form as
\begin{equation}\label{eq-fraLamatrix}
    (-\Delta_h)^{\frac{\alpha}{2}}v = h^{-\alpha}\mathcal{A}^{\alpha}v, 
\end{equation} 
where $\mathcal{A}^{\alpha}$ is a symmetric block Toeplitz matrix with Toeplitz blocks (BTTB), its structure is given by
\begin{equation*}
    \mathcal{A}^{\alpha}=
    \begin{pmatrix}
        \mathcal{A}^{\alpha}_0 & \mathcal{A}^{\alpha}_1 & \mathcal{A}^{\alpha}_2 & \cdots & \mathcal{A}^{\alpha}_{M-2} \\
        \mathcal{A}^{\alpha}_1 & \mathcal{A}^{\alpha}_0 & \mathcal{A}^{\alpha}_1 & \ddots & \mathcal{A}^{\alpha}_{M-3} \\
        \vdots & \mathcal{A}^{\alpha}_1 & \mathcal{A}^{\alpha}_0 & \ddots & \vdots \\
        \mathcal{A}^{\alpha}_{M-3} & \ddots & \ddots & \ddots & \mathcal{A}^{\alpha}_1 \\
        \mathcal{A}^{\alpha}_{M-2} & \mathcal{A}^{\alpha}_{M-3} & \cdots & \cdots & \mathcal{A}^{\alpha}_0
    \end{pmatrix}, \quad
        \mathcal{A}^{\alpha}_j =  
    \begin{pmatrix}
        a_{0,j}^{(\alpha)} & a_{1,j}^{(\alpha)} & a_{2,j}^{(\alpha)} & \cdots & a_{M-2,j}^{(\alpha)} \\
        a_{1,j}^{(\alpha)} & a_{0,j}^{(\alpha)} & a_{1,j}^{(\alpha)} & \ddots & a_{M-3,j}^{(\alpha)} \\
        \vdots & a_{1,j}^{(\alpha)} & a_{0,j}^{(\alpha)} & \ddots & \vdots \\
        a_{M-3,j}^{(\alpha)} & \ddots & \ddots & \ddots & a_{1,j}^{(\alpha)} \\
        a_{M-2,j}^{(\alpha)} & a_{M-3,j}^{(\alpha)} & \cdots & \cdots & a_{0,j}^{(\alpha)}
    \end{pmatrix}.
\end{equation*}

Now, we present some properties of the discrete fractional Laplacian.
\begin{lemma}[Fractional semi-norm equivalence]\label{lem:norm_equivalence}
For $v\in\mathsf C_h^0$, identified with its zero extension,  we have 
\begin{equation}
    \left(\frac{2}{\pi}\right)^\alpha \smnorm{v}_{H_h^{\frac{\alpha}{2}}}^2 \leq \inner{(-\Delta_h)^{\frac{\alpha}{2}} v}{v} \leq \smnorm{v}_{H_h^{\frac{\alpha}{2}}}^2.
\end{equation}

\end{lemma}

\begin{lemma}[A uniform discrete maximum-norm estimate]\label{lem-uniform_maxest}
    For $1<\alpha\leq 2$ and $v\in\mathsf{C}_h^0$, there exists a constant $C_{\alpha}$ independent of $h$, such that 
    \begin{equation}
        \norm{v}_{\infty} \leq C_{\alpha} \left( \norm{v} + \norm{(-\Delta_h)^{\frac{\alpha}{2}} v} \right).
    \end{equation}
\end{lemma}

\begin{proof}
    Let $1<\alpha<2$, and set $s=\alpha/2\in(0,1)$. Using integration by parts and the Gamma integral, for $z\geq0$, one can obtain the identity
    \begin{equation}\label{eq-identity}
        z^s = \frac{s}{\Gamma(1-s)} \int_{0}^{\infty} (1-e^{-rz})r^{-1-s} \rmd r.
    \end{equation}
    Since the generating function can be written as $\varphi(\theta) = 4 - \sum_{\nu=1}^{2}(e^{\bi\theta_{\nu}} + e^{-\bi\theta_{\nu}})$, we can derive
    \begin{equation*}
        e^{-r\varphi(\theta)} = e^{-4r} \prod_{\nu=1}^{2}e^{re^{\bi \theta_{\nu}}} e^{re^{-\bi \theta_{\nu}}} = e^{-4r} \sum_{a,b,c,d \geq 0} \frac{r^{a+b+c+d}}{a!b!c!d!}e^{\bi[(a-b)\theta_1+(c-d)\theta_2]}.
    \end{equation*}
    Take $m=(m_1,m_2) = (a-b,c-d)$ and let $m\cdot\theta=m_1\theta_1+m_2\theta_2$, the above equality can be expressed as $e^{-r\varphi(\theta)} = \sum_{m\in\mathbb{Z}^2}p_{m}(r)e^{\bi m\cdot\theta}$, where
    \begin{equation*}
        p_{m}(r)=e^{-4r}\sum_{\substack{a,b,c,d\geq 0 \\ a-b=m_1,\, c-d=m_2}} \frac{r^{a+b+c+d}}{a!b!c!d!}.
    \end{equation*}
    For $r>0$, we have $p_m(r)=p_{-m}(r)>0$ and  $\sum_{m\in \mathbb{Z}^2} p_{m}(r) = 1$. In particular, we can know  $0\leq 1- p_{\vec{0}}(r) \leq \min\{1,4r\}$  from $p_{\vec{0}}(r) \geq e^{-4r}$, where $\vec{0}=(0,0)$. Then, combining \eqref{eq-compute-intcoeff} and \eqref{eq-identity} and applying Fubini's theorem, we obtain
    \begin{equation}\label{eq-coeident}
        \begin{aligned}
            a_{m}^{(\alpha)} 
            &= \frac{s}{\Gamma(1-s)} \frac{1}{4\pi^2}  \int_{[-\pi,\pi]^2} \left[\int_{0}^{\infty} (1-e^{-r\varphi(\theta)})r^{-1-s} \rmd r \right] e^{-\bi m\cdot \theta} \rmd\theta \\
            &= \frac{s}{\Gamma(1-s)} \int_{0}^{\infty} \left[\frac{1}{(2\pi)^2} \int_{[-\pi,\pi]^2}(1-e^{-r\varphi(\theta)}) e^{-\bi m\cdot\theta} \rmd\theta\right]r^{-1-s} \rmd r \\
            &= \frac{s}{\Gamma(1-s)} \int_{0}^{\infty} \left(\delta_{m \vec{0}} - p_{m}(r)\right)r^{-1-s} \rmd r,
        \end{aligned}
    \end{equation}
    where $\delta_{m\vec0}$ is the Kronecker delta. The bound $|1-e^{-r\varphi(\theta)}|\leq\min\{1,8r\}$ justifies the interchange of integrals in \eqref{eq-coeident}. A simple analysis of \eqref{eq-coeident} shows that
    \begin{equation}\label{eq-coeident-bound}
        a_{m}^{(\alpha)} = a_{-m}^{(\alpha)}\leq 0 \ (m\neq \vec{0}), \qquad a_{\vec{0}}^{(\alpha)} \geq 0, \qquad \sum_{m}a_{m}^{(\alpha)}=0, \qquad \sum_{m\neq \vec{0}}|a_{m}^{(\alpha)}|=a_{\vec{0}}^{(\alpha)}.
    \end{equation}
    By Tonelli's theorem, the above sum identities follow from $\sum_{m\neq\vec0}p_m(r)=1-p_{\vec0}(r)$. The evenness of $\varphi$ in each coordinate also yields $a_{m_1,m_2}^{(\alpha)}=a_{|m_1|,|m_2|}^{(\alpha)}$. At $\alpha=2$, the standard five-point Laplacian has $a_{0,0}^{(2)}=4$, $a_{\pm1,0}^{(2)}=a_{0,\pm1}^{(2)}=-1$, and all other coefficients zero, so \eqref{eq-coeident-bound} remains valid. The rest of the proof applies to $1<\alpha\leq2$.
    
    Now, set $\mathcal{G}_h = \{(i,j)\in\mathbb{Z}^2: 1\leq i,j\leq M-1\}$.
    Recall that the vector $v$ in \eqref{eq-vcolumn} is arranged in column-major order, we define the corresponding index map $\rho: \mathcal{G}_h \to \mathcal{J}_h$ by $\rho(i,j)=i+(j-1)(M-1)$, $(i,j)\in\mathcal{G}_h$, where $\mathcal{J}_h = \rho(\mathcal{G}_h) = \{1,2,\cdots,(M-1)^2\}$. For $p=\rho(i,j)$ and $q=\rho(k,\ell)$, we have 
    \begin{equation*}
        [\mathcal{A}^{\alpha}]_{\rho(i,j),\rho(k,\ell)} = \left[\mathcal{A}_{|j-\ell|}^{\alpha}\right]_{i,k} = a_{|i-k|,|j-\ell|}^{(\alpha)} = a_{i-k,j-\ell}^{(\alpha)},
    \end{equation*} 
    we then derive from \eqref{eq-coeident-bound} that 
    \begin{equation*}
        \sum_{q\in\mathcal{J}_h} [h^{-\alpha}\mathcal{A}^{\alpha}]_{pq} = h^{-\alpha}\sum_{(k,\ell)\in \mathcal{G}_h} a_{i-k,j-\ell}^{(\alpha)} = -h^{-\alpha}\sum_{(k,\ell)\in\mathbb{Z}^2 \backslash\mathcal{G}_h} a_{i-k,j-\ell}^{(\alpha)} \geq 0.
    \end{equation*}
    Let $P(t) = e^{-th^{-\alpha}\mathcal{A}^{\alpha}}$ and $c_h = h^{-\alpha}a_{0,0}^{(\alpha)}$. $P(t)$ can be expressed as 
    \begin{equation*}
        P(t) = e^{-tc_h\mathcal{I}_h}e^{t(c_h\mathcal{I}_h-h^{-\alpha}\mathcal{A}^{\alpha})} = e^{-c_h t}\sum_{n=0}^{\infty}\frac{t^n}{n!}(c_h\mathcal{I}_h-h^{-\alpha}\mathcal{A}^{\alpha})^n,
    \end{equation*}
    where $\mathcal{I}_h$ is the identity matrix of size $(M-1)^2$. By the equation \eqref{eq-coeident-bound}, we know all entries of $c_h\mathcal I_h-h^{-\alpha}\mathcal A^\alpha$ are nonnegative, and its row sums are at most $c_h$. Let $\bm{1}$ be a finite vector in which every component equals $1$. Consequently, $P(t)\bm1\leq\bm1$. By symmetry of $\mathcal A^\alpha$, the column sums of $P(t)$ are also at most $1$.  Introduce the norm $\|w\|_1 = h^2\sum_{p\in\mathcal{J}_h}|w_{p}|$, for $f\in\mathsf{C}_h^0$, we obtain
    \begin{equation}\label{eq-estf1}
        \|P(t)f\|_{1} = h^2\sum_{p\in\mathcal{J}_h} \left|\sum_{q\in\mathcal J_h}[P(t)]_{pq}f_q\right| \leq h^2\sum_{q\in\mathcal{J}_h}|f_q| \sum_{p\in\mathcal{J}_h}[P(t)]_{pq} \leq  h^2\sum_{q\in\mathcal J_h}|f_q| = \|f\|_1.
    \end{equation}

    Denote $E(w) = \inner{(-\Delta_h)^{\frac{\alpha}{2}}w}{w}$, $D_h=[-\pi/h,\pi/h]^2$ and $w^{\ast}(\eta)=h^2\sum_{(i,j)\in\mathcal{G}_h} w_{i,j} e^{-\bi h(\eta_1 i + \eta_2 j)}$. For $R>0$, we first get
    \begin{equation}\label{eq-lowfrees1}
        \frac{1}{4\pi^2} \int_{D_h\cap \{|\eta|\leq R\}}|w^{\ast}(\eta)|^2 \rmd\eta \leq \frac{\|w\|_1^2}{4\pi^2}\mathrm{area}(D_h\cap \{|\eta|\leq R\})  \leq \frac{R^2}{4\pi}\|w\|_1^2,
    \end{equation}
    and then using Lemma \ref{lem:norm_equivalence}, we have
    \begin{equation*}
        \begin{aligned}
        \frac{1}{4\pi^2} \int_{D_h\cap \{|\eta|>R\}} |w^{\ast}(\eta)|^2 \rmd\eta  
        \leq \frac{R^{-\alpha}}{4\pi^2}\int_{D_h\cap\{|\eta|>R\}}|\eta|^\alpha|w^*(\eta)|^2\,d\eta 
        \leq R^{-\alpha}|w|_{H_h^{\frac{\alpha}{2}}}^{\,2}  
        \leq k_\alpha^{-1}R^{-\alpha}E(w),
        \end{aligned}
    \end{equation*}
    where $k_\alpha = (2/\pi)^{\alpha}$. Combining the above two formulas, we obtain
    \begin{equation*}
        \begin{aligned}
            \|w\|^2 
            =& \frac{1}{4\pi^2} \left(\int_{D_h\cap \{|\eta|\leq R\}}+\int_{D_h\cap \{|\eta|>R\}}\right)|w^{\ast}(\eta)|^2 \rmd\eta \leq  \frac{R^2}{4\pi}\|w\|_1^2 + k_\alpha^{-1}R^{-\alpha}E(w),
        \end{aligned}
    \end{equation*}
    For $w\neq0$ and choose $R=(E(w)/\|w\|_1^2)^{1/(\alpha+2)}$, it follows that
    % Raising the resulting bound to the power $(\alpha+2)/2$ gives
    \begin{equation}\label{eq-estf2}
        \|w\|^{\alpha+2} \leq N_{\alpha} E(w)\|w\|_1^{\alpha},\qquad N_{\alpha}=\left(\frac{1}{4\pi}+k_{\alpha}^{-1}\right)^{\frac{\alpha+2}{2}}.
    \end{equation}

    Next, for any $f\in \mathsf{C}_h^0$, we take $w(t)=P(t)f$ and $y(t)=\|w(t)\|^2$. Using \eqref{eq-fraLamatrix}, \eqref{eq-estf1} and \eqref{eq-estf2}, we derive
    \begin{equation*}
        y^{\prime}(t) =  -2\inner{(-\Delta_h)^{\frac{\alpha}{2}}w(t)}{w(t)} = -2E(w(t)) \leq -\frac{2}{N_{\alpha}\|f\|_{1}^{\alpha}}y(t)^{1+\frac{\alpha}{2}}.
    \end{equation*}
    Note that $(y(t)^{-\alpha/2})^{\prime} = (-\alpha/2)y(t)^{-1-\alpha/2}y^{\prime}(t)$. Integrating the inequality for $(y(t)^{-\alpha/2})^{\prime}$, then we have $y(t) \leq \left(N_{\alpha}/\alpha\right)^{2/\alpha}  t^{-2/\alpha}\|f\|_{1}^{2}$, which leads to
    \begin{equation}\label{eq-estf3}
        \|P(t)f\| \leq H_{\alpha} t^{-\frac{1}{\alpha}}\|f\|_{1},\qquad H_{\alpha} = \left(N_{\alpha}/\alpha\right)^{\frac{1}{\alpha}}.
    \end{equation}
    According to the duality identity $\|z\|_\infty=\sup_{\|g\|_1\le1}|\langle z,g\rangle|$ and the self-adjointness of $P(t)$, we obtain
    \begin{equation}\label{eq-maxnorm-es1}
        \|P(t)f\|_\infty=\sup_{\|g\|_1\le1}|\langle f,P(t)g\rangle| \leq \|f\|\sup_{\|g\|_1\le1}\|P(t)g\| \leq H_{\alpha} t^{-\frac{1}{\alpha}} \|f\|,
    \end{equation}
    where the Cauchy-Schwarz inequality and \eqref{eq-estf3} are applied.

    Finally, let $f=(\mathcal{I}_h + (-\Delta_h)^{\frac{\alpha}{2}})v = v + h^{-\alpha}\mathcal{A}^{\alpha}v$, we find that 
    \begin{equation*}
        \frac{\rmd}{\rmd t}\left(e^{-t}P(t)v\right) = -e^{-t}P(t)\left(\mathcal{I}_h + h^{-\alpha}\mathcal{A}^{\alpha}\right)v = -e^{-t}P(t)f,
    \end{equation*}
    which leads to $v = \int_{0}^{\infty} e^{-t}P(t)f \rmd t$. Consequently,
    \begin{equation}\label{eq-resolvemaxes1}
        \begin{aligned}
            \|v\|_{\infty} \leq \int_{0}^{\infty} e^{-t}\|P(t)f\|_{\infty} \rmd t \leq H_{\alpha}\|f\| \int_{0}^{\infty} e^{-t}  t^{-\frac{1}{\alpha}} \rmd t \leq C_{\alpha} \left(\|v\| + \|(-\Delta_h)^{\frac{\alpha}{2}}v\|\right),
        \end{aligned}
    \end{equation}
    where the constant $C_{\alpha}$ is defined as
    \begin{equation*}
        C_{\alpha} = H_{\alpha} \Gamma(1-1/\alpha) = \Gamma(1-1/\alpha)\alpha^{-\frac{1}{\alpha}}\left[\frac{1}{4\pi} + \left(\frac{\pi}{2}\right)^{\alpha}\right]^{\frac{1}{2}+\frac{1}{\alpha}}.
    \end{equation*}
    The integral is finite precisely when $\alpha>1$. This finishes the proof.
\end{proof}

It is easy to conclude from Lemma \ref{lem:norm_equivalence} that $\mathcal{A}^{\alpha}$ is a symmetric positive definite matrix. Consequently, there exists a nonsingular matrix $\Lambda$ such that $\mathcal{A}^{\alpha} = \Lambda^{\top} \Lambda$. This property leads to the following two consequences.
\begin{lemma}\label{lem-fraLa-inneridenty}
    For $v,w\in\mathsf C_h^0$, and let $\mathcal{L}^{\alpha} = h^{-\alpha/2}\Lambda$, we get
    \begin{equation}
        \inner{(-\Delta_h)^{\frac{\alpha}{2}}v}{w} = \inner{v}{(-\Delta_h)^{\frac{\alpha}{2}}w},\quad \inner{(-\Delta_h)^{\frac{\alpha}{2}}v}{w} = \inner{\mathcal{L}^{\alpha}v}{\mathcal{L}^{\alpha}w},
    \end{equation} 
\end{lemma}
\begin{lemma}\label{lem-inner-relipart}
    For any grid function $w^{k}\in \mathsf{C}_h^0$, we obtain 
    \begin{equation}
        \begin{aligned}
            &\repart \inner{(-\Delta_h)^{\frac{\alpha}{2}}\wh{w}^{\,k+\frac{1}{2}}}{\delta_t w^{k+1}} = \frac{1}{2\tau}\left( \norm{\mathcal{L}^{\alpha}w^{k+1}}^2-\norm{\mathcal{L}^{\alpha}w^{k}}^2 \right), \\
                &\impart \inner{(-\Delta_h)^{\frac{\alpha}{2}}\wh{w}^{\,k+\frac{1}{2}}}{\wh{w}^{\, k+\frac{1}{2}}} = 0.
        \end{aligned}
    \end{equation}
\end{lemma}

Define the linear operators $\mathcal{P}_h,\mathcal{Q}_h:\mathsf{C}_h^0 \mapsto \mathsf{C}_h^0$ by
\begin{equation*}
\mathcal{P}_h = \mathcal{I}_h+\frac{\bi \tau}{2}(-\Delta_h)^{\frac{\alpha}{2}},\qquad \mathcal{Q}_h = \mathcal{I}_h-\frac{\bi \tau}{2}(-\Delta_h)^{\frac{\alpha}{2}},    
\end{equation*}
where $\mathcal{I}_h$ is the identity operator on $\mathsf{C}_h^0$. If $\mathcal{Z}_h$ denotes either $\mathcal{P}_h$ or $\mathcal{Q}_h$, Lemma \ref{lem-fraLa-inneridenty} gives
\begin{equation*}
    \repart \inner{\mathcal{Z}_h w}{w} = \norm{w}^2,\quad \forall w \in \mathsf{C}_h^0,
\end{equation*}
which means that $\operatorname{Ker}(\mathcal{Z}_h)=\{0\}$, and thus $\mathcal{Z}_h$ is invertible. Two useful lemmas associated with these linear operators are stated below.
% In the following, we state two crucial lemmas for these linear operators.

\begin{lemma}\label{lem-openorm}
    For the invertible operator $\mathcal{P}_h$, $\mathcal{Q}_h$, and let $\mathcal{B}_h = \mathcal{P}_h^{-1} \mathcal{Q}_h$, they fulfill
    \begin{align}
        &\norm{\mathcal{B}_h w} = \norm{w}, \quad \forall w \in \mathsf{C}_h^0, \label{eq-Bh} \\ 
        &\left(\mathcal{I}_h + \mathcal{B}_h\right)^{-1} \mathcal{P}_h^{-1} =  \frac{1}{2}\mathcal{I}_h. \label{eq-OPide}
    \end{align}  
\end{lemma}
\begin{proof}
    Taking $v=\mathcal B_hw\in\mathsf C_h^0$, we obtain $\mathcal{P}_h v = \mathcal{Q}_h w$ and $v-w = -\frac{\bi\tau}{2}(-\Delta_h)^{\frac{\alpha}{2}}(v+w)$. Using Lemma \ref{lem-fraLa-inneridenty}, we have $\repart \inner{v-w}{v+w}=0$, and hence $\norm{v} = \norm{w}$. Thus, \eqref{eq-Bh} is confirmed.  From the definition of these operators, we know that 
    \begin{equation*}
        \mathcal{P}_h(\mathcal{I}_h+ \mathcal{B}_h) = \mathcal{P}_h + \mathcal{Q}_h = 2\mathcal{I}_h.
    \end{equation*}
    Taking the inverse of both sides of the above equation, then \eqref{eq-OPide} is established.
\end{proof}

\begin{lemma}\label{lem-opeequality}
    Let $\{w^k\}_{k=1}^{K}$ be a sequence in $\mathsf{C}_h^0$, and it satisfies
    \begin{equation}\label{eq-opeequcond}
        w^{k+1} = \left(\mathcal{B}_h - \mathcal{I}_h\right)w^k + \mathcal{B}_h w^{k-1} + \mathcal{P}_h^{-1}\mathsf{\Gamma}^k,\quad 2\leq k\leq K-1.
    \end{equation}
    Then, for $2\leq k\leq K-1$, we have
    \begin{equation}\label{eq-opeequres}
        \begin{aligned}
            w^{p} =& \frac{1}{2}\left[(-1)^{p}\mathcal{I}_h + \mathcal{B}_h^{p-1}\right]\mathcal{P}_h w^2 + \frac{1}{2}\left[(-1)^{p-1}\mathcal{B}_h + \mathcal{B}_h^{p-1}\right]\mathcal{P}_h w^1 \\ 
            &+ \frac{1}{2}\sum_{l=2}^{k}\left[(-1)^{p-l+1}\mathcal{I}_h + \mathcal{B}_h^{p-l} \right]\mathsf{\Gamma}^{l}, \qquad p=k,k+1.
        \end{aligned}
    \end{equation} 
\end{lemma}
\begin{proof}
    It is easy to know that \eqref{eq-opeequcond} can be expressed in a vector form
    \begin{equation*}
        \begin{bmatrix} 
            w^{k+1} \\
            w^{k} 
        \end{bmatrix}
        =
        S
        \begin{bmatrix}
            w^{k} \\
            w^{k-1} 
        \end{bmatrix}
        +
        \begin{bmatrix}
            \mathcal{P}_h^{-1}\mathsf{\Gamma}^k \\
            0
        \end{bmatrix}, \qquad S=
        \begin{bmatrix}
            \mathcal{B}_h - \mathcal{I}_h & \mathcal{B}_h \\
            \mathcal{I}_h & 0
        \end{bmatrix}.
    \end{equation*}
    By a simple induction argument, we obtain
    \begin{equation}\label{eq-opefm1}
           \begin{bmatrix} 
            w^{k+1} \\
            w^{k} 
        \end{bmatrix}
        =
        S^{k-1}
        \begin{bmatrix}
            w^{2} \\
            w^{1} 
        \end{bmatrix}
        +\sum_{l=2}^{k}S^{k-l}
        \begin{bmatrix}
            \mathcal{P}_h^{-1}\mathsf{\Gamma}^l \\
            0
        \end{bmatrix}.
    \end{equation}
    According to the Part 9 of the proof of Theorem 3.1 in \cite{zouraris2023error} and \eqref{eq-OPide}, we get
    \begin{equation}\label{eq-opefm2}
        S^{j}=\frac{1}{2}
        \begin{bmatrix}
            [(-1)^{j}\mathcal{I}_h + \mathcal{B}_h^{j+1}]\mathcal{P}_h \hspace{2pt} & \left[(-1)^{j+1}\mathcal{B}_h + \mathcal{B}_h^{j+1}\right]\mathcal{P}_h \\
           [(-1)^{j+1}\mathcal{I}_h + \mathcal{B}_h^{j}]\mathcal{P}_h \hspace{2pt} & \left[(-1)^{j}\mathcal{B}_h + \mathcal{B}_h^{j}\right]\mathcal{P}_h            
        \end{bmatrix}
        ,\qquad j=1,2,\cdots.
    \end{equation}
    By plugging \eqref{eq-opefm2} into \eqref{eq-opefm1}, we derive \eqref{eq-opeequres}.
    % \tcr{The proof is omitted for simplicity, interested readers can consult Part 9 of Theorem 3.1 in \cite{zouraris2023error}.}
\end{proof}

\section{The relaxation Crank-Nicolson difference scheme}\label{sec-scheme}

In this section, we first establish the RCND scheme for the SFNS equation, and then derive the truncation error through the consistency equation. Finally, the unique solvability and discrete conservation laws are discussed.

\subsection{Fully discrete scheme}
Introduce the relaxation variable $\psi:=g(|u|^2)$. Then, the SFNS equation \eqref{eq-SFNS} with the boundary condition \eqref{eq-bondcd} can be reformulated as 

\begin{equation}\label{eq-relaxNS}
    \begin{cases}
        \bi u_t - (-\Delta)^{\frac{\alpha}{2}} u + \psi u = 0, \quad &(x,y) \in \Omega, \quad 0<t<T,\\
        \psi = g(|u|^2).
    \end{cases}
\end{equation}

Next, we perform full discretization on the equation \eqref{eq-relaxNS}.
We employ the Crank-Nicolson time stepping for $u$ at integer time levels, and a staggered relaxation recurrence scheme for $\psi$ at half-integer levels. In the spatial direction, a second-order fractional central difference scheme in \eqref{eq-fulldiscrete-fraLa} is used to approximate the fractional Laplacian; i.e., we replace $(-\Delta)^{\frac{\alpha}{2}}$ with $(-\Delta_h)^{\frac{\alpha}{2}}$. To preserve the global second-order temporal accuracy of the relaxation variable, an initialization procedure is employed to obtain the solution $U^{1/2}$ at the half time level $t_{1/2}$. The RCND scheme is detailed in the following three steps:

\textbf{Step 1.} Given the initial values $U^0:=u^0$ and $\Psi^0:=g(|u^0|^2)$, solve for $U^{1/2}\in \mathsf{C}_h^0$ satisfying
\begin{equation}
    \bi\frac{U^{\frac{1}{2}}-U^0}{\tau/2} - (-\Delta_h)^{\frac{\alpha}{2}}\frac{U^{\frac{1}{2}}+U^0}{2} + \Psi^0\frac{U^{\frac{1}{2}}+U^0}{2}=0. \label{eq-numsma} 
\end{equation} 

\textbf{Step 2.} Let $\Psi^{1/2} := g(|U^{1/2}|^{2}) \in \mathsf{R}_h$, and compute $U^1 \in \mathsf{C}_h^0$ from 
\begin{equation}
    \bi\delta_t U^1 - (-\Delta_h)^{\frac{\alpha}{2}} \wh{U}^{\frac{1}{2}} + \Psi^{\frac{1}{2}}\wh{U}^{\frac{1}{2}}=0. \label{eq-numsmb}
\end{equation}

\textbf{Step 3.} For $k=1,\cdots,K-1$, update $\Psi^{k+1/2} \in \mathsf{R}_h$ by
\begin{equation}\label{eq-numsmc}
    \Psi^{k+\frac{1}{2}}:=2g(| U^{k}|^2)-\Psi^{k-\frac{1}{2}},
\end{equation}
and then evaluate $U^{k+1} \in \mathsf{C}_h^0 $ by solving the equation
\begin{equation}
    \bi\delta_t U^{k+1} - (-\Delta_h)^{\frac{\alpha}{2}} \wh{U}^{k+\frac{1}{2}} + \Psi^{k+\frac{1}{2}}\wh{U}^{k+\frac{1}{2}}=0. \label{eq-numsmd} 
\end{equation}

\begin{remark}
    If Step 1 is omitted and $\Psi^{1/2}$ is initialized by $g(|u^0|^2)$, then, in the cubic case, \eqref{eq-numsmb}--\eqref{eq-numsmd} reduces to the Besse relaxation scheme \cite{besse2004relaxation}. Although this scheme can achieve second-order accurate at integer time levels, the relaxation variable at half-integer time levels is only first-order accurate \cite{zouraris2023error}. The designed scheme \eqref{eq-numsma}--\eqref{eq-numsmd} can overcome this limitation.
\end{remark}

\subsection{Consistency error analysis}

Let $u^{1/4}$, $u^k$, $u^{k+1/2} \in \mathsf{C}_h^0$ denote the exact solution at $t_{1/4}$, $t_k$ and $t_{k+1/2}$, respectively. 
Similarly, let $\psi^{1/4}, \psi^{k}, \psi^{k+1/2} \in \mathsf{R}_h$ denote the exact relaxation variables at the same time levels.
Then, substituting them into the discrete system \eqref{eq-numsma}--\eqref{eq-numsmd}, we obtain the following consistency equations:

\begin{itemize}
    \item At $t=t_{1/4}$, %we have
        \begin{equation}\label{eq-consiseqa}  
            \bi \frac{u^{\frac{1}{2}}-u^{0}}{\tau/2} - (-\Delta_h)^{\frac{\alpha}{2}} \frac{u^{\frac{1}{2}}+u^{0}}{2} + \psi^{0} \frac{u^{\frac{1}{2}}+u^{0}}{2} = \mathtt{r}_u^{\frac{1}{4}}. 
        \end{equation}
    \item At $t=t_{1/2}$, with $\psi^{1/2} = g(|u^{1/2}|^2)$, %we have
        \begin{equation}\label{eq-consiseqb} 
            \bi \delta_t u^{1} - (-\Delta_h)^{\frac{\alpha}{2}} \wh{u}^{\frac{1}{2}} + \psi^{\frac{1}{2}}\wh{u}^{\frac{1}{2}}  = \mathtt{r}_u^{\frac{1}{2}}. 
        \end{equation}
    \item  At $t=t_{k}$, % we get 
        \begin{equation}\label{eq-consiseqc}
            \psi^{k+\frac{1}{2}} = 2g(|u^k|^2) - \psi^{k-\frac{1}{2}} + \mathtt{r}_{\psi}^{k}, \qquad 1\leq k\leq K-1,      
        \end{equation}
        and at $t=t_{k+1/2}$, %we obtain
        \begin{equation}\label{eq-consiseqd}
            \bi \delta_t u^{k+1} - (-\Delta_h)^{\frac{\alpha}{2}} \wh{u}^{\, k+\frac{1}{2}} + \psi^{k+\frac{1}{2}}\wh{u}^{\, k+\frac{1}{2}}  = \mathtt{r}_u^{k+\frac{1}{2}}, \qquad 1\leq k\leq K-1.    
        \end{equation} 
\end{itemize}
Among the above equations, $\mathtt{r}_u^{1/4}$, $\mathtt{r}_{\psi}^{k}$ and $\mathtt{r}_u^{k+1/2}$ are the truncation errors at their corresponding time levels, whose expressions are given as follows:
% Introduce a indicator variable $\sigma = \{u,v\}$, $\varsigma = \{\psi,\phi\}$ and let  $\Upsilon_u = \psi + \rho\phi$, $\Upsilon_v = \phi + \rho\psi$. Then, truncation error terms are given by
\begin{equation*}
    \begin{aligned}
        &\mathtt{r}_u^{\frac{1}{4}} := \bi \left(\frac{u^{\frac{1}{2}} - u^0}{\tau/2} - \partial_t u^{\frac{1}{4}} \right) - \left[ (-\Delta_h)^{\frac{\alpha}{2}} \frac{u^{\frac{1}{2}} + u^0}{2}  - (-\Delta)^{\frac{\alpha}{2}} u^{\frac{1}{4}} \right] \\
        &\hspace{1em}  + \left(\psi^0 \frac{u^{\frac{1}{2}} + u^0}{2} - \psi^{\frac{1}{4}} u^{\frac{1}{4}} \right) , \\
        &\mathtt{r}_u^{k+\frac{1}{2}} := \bi \bigl(\delta_t u^{k+1} - \partial_t u^{k+\frac{1}{2}}\bigr) - \left[(-\Delta_h)^{\frac{\alpha}{2}}\wh{u}^{\, k+\frac{1}{2}} - (-\Delta)^{\frac{\alpha}{2}} u^{k+\frac{1}{2}} \right]  \\
        &\begin{array}{@{}l @{\qquad}l@{}}
        \hspace{1.3em} + \left[ \psi^{k+\frac{1}{2}}\bigl(\wh{u}^{\, k+\frac{1}{2}} - u^{k+\frac{1}{2}} \bigr) \right] ,  &  0 \leq k \leq K-1, \\
        \mathtt{r}_\psi^k := \psi^{k+\frac{1}{2}} + \psi^{k-\frac{1}{2}} - 2\psi^k , &  1 \leq k \leq K-1.
        \end{array}
    \end{aligned} 
\end{equation*}

The following lemma gives estimates for the above truncation errors. 
% under suitable regularity assumptions.

\begin{lemma}\label{lem-truncerror}
    Assume $u\in C^{4}\left([0,T],\mathcal{W}^{2+\alpha}(\mathbb R^2)\right)$ and $g\in C^3(\mathbb{R};\mathbb{R})$. 
    \begin{itemize}
        \item The above truncation errors satisfy
        \begin{equation} \label{eq-truncerr1}
            \bigl\| \mathtt{r}_u^{\frac{1}{4}}\bigr\|  \leq \widetilde{C}\left(\tau+h^2\right), \quad
            \max_{0 \leq k\leq K-1}\bigl\|\mathtt{r}_u^{k+\frac{1}{2}}\bigr\|  \leq \widetilde{C}\left(\tau^2+h^2\right), \quad
            \max_{1\leq k\leq K-1} \left\|\mathtt{r}_\psi^k \right\|+\left\|\mathtt{r}_\psi^k\right\|_\infty \leq \widetilde{C}\tau^2.
        \end{equation} 
        \item Their temporal differences satisfy
        \begin{equation}\label{eq-truncerr2}
            \max_{1\leq k\leq K-1}\bigl\|\mathtt{r}_u^{k+\frac{1}{2}}-\mathtt{r}_u^{k-\frac{1}{2}}\bigr\| \leq \widetilde{C}\tau\left(\tau^2+h^2\right), \quad 
            \max_{2\leq k\leq K-1}\bigl\|\mathtt{r}_\psi^k-\mathtt{r}_\psi^{k-1}\bigr\| \leq \widetilde{C}\tau^3,
        \end{equation}
    \end{itemize}
    where $\widetilde{C}$ is a positive constant independent of $\tau$ and $h$.
\end{lemma}

\begin{proof}
    Combining the Taylor expansion formula with Lemma \ref{lem:discrete_laplacian}, we get 
    \begin{equation*}
        \begin{aligned}
            (\mathtt{r}_u)^{\frac{1}{4}}_{i,j} =& \bi\frac{\tau^2}{16}\int_{-\frac{1}{4}}^{\frac{1}{4}}\left(1-4|s|\right)^2 \partial_t^3 u\bigl(x_i,y_j,t_{\frac{1}{4}}+\tau s \bigr) \, \rmd s - \mathcal{O}(h^2) \\
            &-\frac{\tau^2}{8}\int_{-\frac{1}{4}}^{\frac{1}{4}} \left(1 - 4|s|\right)  \left[(-\Delta)^{\frac{\alpha}{2}}\partial_t^2 u\bigl(x_i,y_j,t_{\frac{1}{4}}+\tau s \bigr) \right] \rmd s \\
            & +\psi^0_{i,j} \frac{\tau^2}{8}\int_{-\frac{1}{4}}^{\frac{1}{4}} \left(1 - 4|s|\right) \,  \partial_t^2 u\bigl(x_i,y_j,t_{\frac{1}{4}}+\tau s\bigr) \, \rmd s -u_{i,j}^{\frac{1}{4}}\int_0^{\tau/4}\partial_t\psi(x_i,y_j,s)\,\rmd s
        \end{aligned}
    \end{equation*}
    for $1\leq i,j\leq M-1$, and
    \begin{equation*}
        \begin{aligned}
            (\mathtt{r}_u)^{k+\frac{1}{2}}_{i,j} =& \bi\frac{\tau^2}{8}\int_{-\frac{1}{2}}^{\frac{1}{2}} \left(1-2|s|\right)^2 \partial_t^3 u\bigl(x_i,y_j,t_{k+\frac{1}{2}}+\tau s\bigr) \, \rmd s - \mathcal{O}(h^2) \\
            &-\frac{\tau^2}{4}\int_{-\frac{1}{2}}^{\frac{1}{2}} \left(1 - 2|s|\right)  \left[(-\Delta)^{\frac{\alpha}{2}}\partial_t^2 u\bigl(x_i,y_j,t_{k+\frac{1}{2}}+\tau s\bigr) \right] \rmd s  \\
            &+\psi^{k+\frac{1}{2}}_{i,j}\frac{\tau^2}{4}\int_{-\frac{1}{2}}^{\frac{1}{2}} \left(1 - 2|s|\right) \, \partial_t^2 u\bigl(x_i,y_j,t_{k+\frac{1}{2}}+\tau s \bigr) \, \rmd s 
        \end{aligned}
    \end{equation*}
    for $1\leq i,j\leq M-1$ and $0\leq k\leq K-1$, which obviously produce the first two estimates of \eqref{eq-truncerr1}. By using the standard Taylor formula again, we obtain
    \begin{equation*}
        (\mathtt{r}_\psi)^k_{i,j} = \frac{\tau^2}{2}\int_{-\frac{1}{2}}^{\frac{1}{2}} \left(1 - 2|s|\right) \, \partial_t^2 \psi \bigl(x_i,y_j,t_{k}+\tau s \bigr) \, \rmd s, \qquad 1 \leq k \leq K-1,
    \end{equation*}
    which yields the last estimate in \eqref{eq-truncerr1}. Similarly, we can also derive the difference estimates \eqref{eq-truncerr2}.
\end{proof}

\subsection{Solvability and discrete conservation laws}

% In this subsection, 
Next, we analyze the unique solvability of the proposed RCND scheme \eqref{eq-numsma}--\eqref{eq-numsmd}.

\begin{theorem}\label{thm-uniquesolve}
    The RCND scheme \eqref{eq-numsma}--\eqref{eq-numsmd} is uniquely solvable.
\end{theorem}

\begin{proof}
Consider the homogeneous system of the equation \eqref{eq-numsma}, that is
\begin{equation}\label{eq-homsys1}
    \bi U^{\frac{1}{2}} - \frac{\tau}{4}(-\Delta_h)^{\frac{\alpha}{2}}U^{\frac{1}{2}} + \frac{\tau}{4} \Psi^0 U^{\frac{1}{2}} = 0.
\end{equation}
Taking the inner product of \eqref{eq-homsys1} with $U^{1/2}$ and using Lemma \ref{lem-fraLa-inneridenty}, we have 
\begin{equation}\label{eq-homsys_norm1}
    \bi \bigl\| U^{\frac{1}{2}}\bigr\|^2 - \frac{\tau}{4} \bigl\| \mathcal{L}^{\alpha}U^{\frac{1}{2}} \bigr\|^2
    + \frac{\tau}{4} \inner{\Psi^0 U^{\frac{1}{2}}}{U^{\frac{1}{2}}}  = 0. 
\end{equation}
Then, taking the imaginary part  of \eqref{eq-homsys_norm1} gives $\| U^{1/2} \|^2 = 0$. Therefore, the equation \eqref{eq-homsys1} has only the zero solution, which means that the system \eqref{eq-numsma} admits a unique solution.

The same argument applies to Step 2, since $\Psi^{1/2}$ is real. Now assume that $\{U^k\}_{k=1}^n$ are known, and the recurrence relation \eqref{eq-numsmc} determines $\{\Psi^{k+1/2}\}_{k=1}^{n}$. Consider the homogeneous system of the equation \eqref{eq-numsmd} for $k=n$, 
\begin{equation}\label{eq-homsys2}
    \bi U^{n+1} - \frac{\tau}{2}(-\Delta_h)^{\frac{\alpha}{2}} U^{n+1} + \frac{\tau}{2}\Psi^{n+\frac{1}{2}} U^{n+1} = 0.
\end{equation}
Taking the inner product of \eqref{eq-homsys2} with $U^{n+1}$ and utilizing Lemma \ref{lem-fraLa-inneridenty}, we get
\begin{equation}\label{eq-homsys_norm2}
    \bi \left\|U^{n+1}\right\|^2 - \frac{\tau}{2}\left\|\mathcal{L}^{\alpha}U^{n+1} \right\|^2 + \frac{\tau}{2} \inner{\Psi^{n+\frac{1}{2}}U^{n+1}}{U^{n+1}} = 0,
\end{equation}
where $\Psi^{n+\frac{1}{2}}$ is a real vector. 
It is easy to obtain $\|U^{n+1}\|^2=0$ by taking the imaginary part of \eqref{eq-homsys_norm2}, so the equation \eqref{eq-numsmd} admits a unique solution.
% By taking the imaginary part of the \eqref{eq-homsys_norm2}, we obtain $\|U^{n+1}\|^2=0$, so the equation \eqref{eq-numsmd} exists a unique solution.

Consequently, the RCND scheme \eqref{eq-numsma}--\eqref{eq-numsmd} is uniquely solvable at each time level. 
\end{proof}

Note that \eqref{eq-numsmc} defines a linear relaxation recurrence, which does not affect the mass conservation of the RCND scheme. However, due to the linear structure of \eqref{eq-numsmc}, the energy conservation is limited to the cubic case. The relevant conservation properties are presented as follows.

\begin{theorem}\label{thm-conservation}
    The RCND scheme \eqref{eq-numsma}--\eqref{eq-numsmd} satisfies the following discrete conservation laws.
    \begin{itemize}
        \item For a general real-valued function $g$, let $\mathtt{M}^k := \|U^k\|^2$, we have the discrete mass conservation: 
        \begin{equation}
            \mathtt{M}^k  =  \mathtt{M}^0, \qquad 0\leq k\leq K.
        \end{equation}
        \item When $g(s)=\beta s$ with $\beta\in\mathbb{R}\backslash \{0\}$, we obtain the discrete mass and energy conservation: 
        \begin{equation}
            \mathtt{M}^k=\mathtt{R}^k=\mathtt{M}^0, \qquad \mathtt{E}^k = \mathtt{E}^0,\qquad 0\leq k\leq K,
        \end{equation}
        where 
        \begin{equation*}
            \mathtt{R}^k:= \frac{1}{2\beta}\left\langle\Psi^{k+\frac{1}{2}}+\Psi^{k-\frac{1}{2}},1\right\rangle,\qquad \mathtt{E}^k := \norm{\mathcal{L}^{\alpha}U^k}^2 - \frac{1}{2\beta} \left\langle\Psi^{k-\frac{1}{2}},\Psi^{k+\frac{1}{2}} \right\rangle.
        \end{equation*}
        Here, $\Psi^{-1/2}:= 2\beta|U^{0}|^2 - \beta|U^{1/2}|^2$ and $\Psi^{K+1/2}:=2\beta|U^K|^2-\Psi^{K-1/2}$.
    \end{itemize}
\end{theorem}

\begin{proof}
The proof is divided into two steps.
\begin{enumerate}[label=(\Roman*)]
    \item Taking the inner product of \eqref{eq-numsmd} with $\wh{U}^{k+\frac{1}{2}}$ and using Lemma \ref{lem-fraLa-inneridenty}, we have
         \begin{equation*}
            \bi \inner{\delta_t U^{k+1}}{\wh{U}^{k+\frac{1}{2}}}  - \bigl\| \mathcal{L}^{\alpha}\wh{U}^{k+\frac{1}{2}} \bigr\|^2 + \inner{\Psi^{k+\frac{1}{2}} \wh{U}^{k+\frac{1}{2}}}{\wh{U}^{k+\frac{1}{2}}}= 0,
        \end{equation*}
        its imaginary part is
        \begin{equation*}
            \frac{1}{2\tau}\left( \left\|U^{k+1} \right\|^2 - \left\|U^{k} \right\|^2 \right) = 0,
        \end{equation*}
        which implies that $\mathtt{M}^k = \mathtt{M}^{k-1} = \cdots =\mathtt{M}^0$ holds.  
    \item In the case of  $g(s)=\beta s$ ($\beta\neq 0$),  it can be seen from \eqref{eq-numsmc} that 
        \begin{equation*}
            \left|U^{k}\right|^2 = \frac{\Psi^{k+\frac{1}{2}} + \Psi^{k-\frac{1}{2}}}{2\beta}, \qquad 1\leq  k\leq K.
        \end{equation*}
        According to the definition of the discrete inner product, we get
        \begin{equation*}
            \mathtt{M}^k = \left\| U^{k} \right\|^2 = \left\langle\left| U^{k} \right|^2,1 \right\rangle = \frac{1}{2\beta} \left\langle\Psi^{k+\frac{1}{2}}+\Psi^{k-\frac{1}{2}}, 1\right\rangle = \mathtt{R}^k, \qquad 1\leq  k\leq K.
        \end{equation*}
        In particular, for $k=0$, it follows that 
        \begin{equation*}
            \mathtt{M}^0 = \left\| U^{0} \right\|^2 = \left\langle\left| U^{0} \right|^2 ,1 \right\rangle =  \frac{1}{2\beta} \left\langle \beta \bigl| U^{\frac{1}{2}} \bigr|^2 + \bigl(2\beta\left| U^{0} \right|^2 - \beta\bigl| U^{\frac{1}{2}} \bigr|^2 \bigr) , 1 \right\rangle 
            = \frac{1}{2\beta} \left\langle \Psi^{\frac{1}{2}}+\Psi^{-\frac{1}{2}} , 1\right\rangle = \mathtt{R}^0.
        \end{equation*}

        Then, taking the inner product of \eqref{eq-numsmd} with $\delta_t U^{k+1}$ and taking its real part, we have
        \begin{equation*}
            \repart \left\langle (-\Delta_h)^{\frac{\alpha}{2}} \wh{U}^{k+\frac{1}{2}} ,\delta_t U^{k+1} \right\rangle
            =\repart\left\langle\Psi^{k+\frac{1}{2}}\wh{U}^{k+\frac{1}{2}},\delta_tU^{k+1}\right\rangle, \qquad 0\leq k\leq K-1.
        \end{equation*}
        % where $0\leq k\leq K-1$. 
        By using Lemmas \ref{lem-fraLa-inneridenty} and \ref{lem-inner-relipart}, the above equality can be simplified as
        \begin{equation}\label{eq-enerymid1}
            \begin{aligned}
                 \norm{\mathcal{L}^{\alpha}U^{k+1}}^2 -\norm{\mathcal{L}^{\alpha}U^{k}}^2  
                = \left\langle \Psi^{k+\frac{1}{2}} , \left| U^{k+1} \right|^2 - \left| U^{k} \right|^2 \right\rangle.
            \end{aligned}
        \end{equation}
        Furthermore, the relaxation recurrence \eqref{eq-numsmc} gives 
        \begin{equation*}
            \left| U^{k+1} \right|^2 - \left| U^{k} \right|^2  = \frac{\Psi^{k+\frac{3}{2}} - \Psi^{k-\frac{1}{2}}}{2\beta},\quad 1\leq k \leq K-1.
        \end{equation*}
        Thus, for $1\leq k \leq K-1$, we obtain
        \begin{equation}\label{eq-enerymid2}
            \left\langle \Psi^{k+\frac{1}{2}} , \left| U^{k+1} \right|^2 - \left| U^{k} \right|^2 \right\rangle  
            = \frac{1}{2\beta} \left[\left\langle \Psi^{k+\frac{3}{2}}, \Psi^{k+\frac{1}{2}} \right\rangle - \left\langle \Psi^{k+\frac{1}{2}}, \Psi^{k-\frac{1}{2}}\right\rangle \right] . 
        \end{equation}
        By substituting \eqref{eq-enerymid2} into the right-hand side of \eqref{eq-enerymid1}, $\mathtt{E}^{k+1} = \mathtt{E}^{k} = \cdots = \mathtt{E}^{1}$ holds.
        In particular, because of 
        \begin{equation*}
            2\beta\left( \left| U^{1} \right|^2 - \left| U^{0} \right|^2 \right) = \beta\left( 2\left| U^{1} \right|^2 - \bigl| U^{\frac{1}{2}} \bigr|^2 \right) - \beta\left( 2\left| U^{0} \right|^2 - \bigl| U^{\frac{1}{2}} \bigr|^2 \right) = \Psi^{\frac{3}{2}} - \Psi^{-\frac{1}{2}}, 
        \end{equation*}
        so the substitution leading to \eqref{eq-enerymid2} is also valid at $k=0$, and hence $\mathtt E^1=\mathtt E^0$.
\end{enumerate}
This finishes the proof.
\end{proof}

\section{Unconditional convergence analysis}\label{sec-convergence}
In this section, our goal is to establish the unconditional optimal-order convergence of the RCND scheme \eqref{eq-numsma}--\eqref{eq-numsmd}. 
To avoid dealing with the discrete $L^\infty$-norm bound for the numerical solution, we introduce two cut-off coefficients and design a modified error system. 
Next, we first demonstrate the closed estimates for the associated auxiliary errors, and then show that the cut-off condition is inactive. 
The resulting estimates apply to the original error system without any mesh-ratio restriction.

\subsection{The modified error system}
Set the error functions as 
\begin{equation*}
    \begin{aligned}
      &e_u^{\frac{1}{2}}=u^{\frac{1}{2}}-U^{\frac{1}{2}}, &&\quad  e_u^k=u^k-U^k, &&\quad e_\psi^{k+\frac{1}{2}}=\psi^{k+\frac{1}{2}}-\Psi^{k+\frac{1}{2}}. 
    \end{aligned}
\end{equation*}
Clearly, $e_u^0=0$. By subtracting \eqref{eq-numsma}--\eqref{eq-numsmd} from \eqref{eq-consiseqa}--\eqref{eq-consiseqd}, we obtain the following error equations.
For the initial step, the errors satisfy
\begin{equation}\label{eq-orierreq1}
    \bi e_u^{\frac{1}{2}}-\frac{\tau}{4}(-\Delta_h)^{\frac{\alpha}{2}} e_u^{\frac{1}{2}}+\frac{\tau}{4}\Psi^0 e_u^{\frac{1}{2}}=\frac{\tau}{2}\mathtt{r}_u^{\frac{1}{4}},  
\end{equation}
and
\begin{equation}\label{eq-orierreq3}
    e_\psi^{\frac{1}{2}}=g\left(\bigl| u^{\frac{1}{2}} \bigr|^2 \right)-g\left(\bigl| U^{\frac{1}{2}} \bigr|^2\right).
\end{equation} 
For $0\leq k\leq K-1$, the primary variable errors satisfy
\begin{equation}\label{eq-orierreq4}
    \bi\delta_t e_u^{k+1}-(-\Delta_h)^{\frac{\alpha}{2}} \wh{e}_u^{\,k+\frac{1}{2}}=\mathrm{F}_{u,1}^{k+\frac{1}{2}}+\mathrm{F}_{u,2}^{k+\frac{1}{2}}+\mathrm{F}_{u,3}^{k+\frac{1}{2}}, 
\end{equation}
where 
\begin{equation*}
    \mathrm{F}_{u,1}^{k+\frac{1}{2}}=-e_\psi^{k+\frac{1}{2}}\wh{U}^{k+\frac{1}{2}}, \qquad 
    \mathrm{F}_{u,2}^{k+\frac{1}{2}}=-\psi^{k+\frac{1}{2}}\wh{e}_u^{\,k+\frac{1}{2}}, \qquad 
    \mathrm{F}_{u,3}^{k+\frac{1}{2}}= \mathtt{r}_{u}^{k+\frac{1}{2}}.
\end{equation*}
For $1\leq k\leq K-1$, the relaxation variable errors satisfy
\begin{equation}\label{eq-orierreq6}
    e_\psi^{k+\frac{1}{2}}+e_\psi^{k-\frac{1}{2}}
    =2\left[g\left(\bigl| u^k \bigr|^2 \right)-g\left(\bigl| U^k \bigr|^2 \right)\right]+\mathtt{r}_\psi^k.
\end{equation}

In what follows, we assume that the exact solutions of SFNS equation \eqref{eq-SFNS} are uniformly bounded; i.e., there exists a  constant $\mathcal{N}>0$ such that  
\begin{equation}\label{eq-exsolbound}
    % \max_{\substack{(x,y)\in \Omega, \\ t\in[0,T]  }} 
    \max_{(x,y)\in \Omega, \, t\in[0,T]  }
    \left\{\left| \varpi(x,y,t) \right|, \left| \partial_t\varpi(x,y,t) \right|\right\} \leq \mathcal{N} , \qquad \varpi= u,\psi.
\end{equation}
Throughout the analysis, $C$ denotes a generic positive constant whose value may change across different
instances but is independent of $\tau$ and $h$.  

Now, we first present an estimate for the error function \eqref{eq-orierreq3}.

\begin{lemma}\label{lem-halfrexest}
Under the assumptions of Lemma \ref{lem-truncerror}, and let $U^{1/2}\in\mathsf{C}_h^0$ be a solution to equation \eqref{eq-numsma}. Then, the initial relaxation error $e_\psi^{1/2}\in \mathsf{R}_h^0$ satisfies
\begin{equation}\label{eq-miderrest}
    \bigl\| e_{\psi}^{\frac{1}{2}} \bigr\| \leq \widetilde{C}_0\left(\tau^2 + h^2 \right) \qquad \text{and} 
    \qquad \bigl\| e_{\psi}^{\frac{1}{2}} \bigr\|_{\infty} \leq \widetilde{C}_0\left(\tau + h^2 \right),
\end{equation}
where $\widetilde{C}_0$ is a positive constant independent of $\tau$ and $h$.
\end{lemma}

\begin{proof}
Taking the inner product of \eqref{eq-orierreq1} with $e_u^{1/2}$, and choosing its imaginary part, we obtain
\begin{equation*}
    \bigl\| e_u^{\frac{1}{2}} \bigr\|^2 =   \frac{\tau}{2}\impart \left\langle\mathtt{r}_u^{\frac{1}{4}}, e_u^{\frac{1}{2}} \right\rangle.
\end{equation*}
By employing Lemma \ref{lem-truncerror} in tandem with Cauchy-Schwarz inequality, we have
\begin{equation}\label{eq-mides1}
    \bigl\| e_u^{\frac{1}{2}} \bigr\| \leq \frac{\tau}{2} \bigl\| \mathtt{r}_u^{\frac{1}{4}} \bigr\| \leq \widetilde{C}\tau\left(\tau + h^2\right).
\end{equation}
Afterwards, it can be concluded from \eqref{eq-orierreq1} and the uniform bound for $\Psi^0=g(|u^0|^2)$ that
\begin{equation}\label{eq-mides2}
    \left\|(-\Delta_h)^{\frac{\alpha}{2}}e_u^{\frac{1}{2}} \right\| \leq C \left(\tau^{-1} \bigl\| e_u^{\frac{1}{2}} \bigr\| + \bigl\|e_u^{\frac{1}{2}} \bigr\| + \bigl\| \mathtt{r}_u^{\frac{1}{4}} \bigr\|\right) \leq C \left(  \tau + h^2 \right),
\end{equation}
and combining \eqref{eq-mides1}--\eqref{eq-mides2} with Lemma \ref{lem-uniform_maxest}, we have
\begin{equation*}
    \bigl\| e_u^{\frac{1}{2}} \bigr\|_\infty \leq C_{\alpha}  \left(\bigl\| e_u^{\frac{1}{2}} \bigr\| + \left\|(-\Delta_h)^{\frac{\alpha}{2}}e_u^{\frac{1}{2}} \right\| \right) \leq C \left(  \tau + h^2 \right).
\end{equation*}
% which implies that numerical solution $U^{\frac{1}{2}}$ meets the uniform boundedness. 
Thus, for sufficiently small $\tau$ and $h$, we get
\begin{equation}\label{eq-midnumbound}
    \bigl\| U^{\frac{1}{2}} \bigr\|_{\infty} \leq  \bigl\| u^{\frac{1}{2}} \bigr\|_\infty + \bigl\| e_u^{\frac{1}{2}} \bigr\|_\infty \leq \mathcal{N} + 1.
\end{equation}
Finally, it follows from \eqref{eq-orierreq3}, \eqref{eq-mides1}, \eqref{eq-midnumbound} and Lemmas \ref{lem-gfunLip} (with $\sigma = (\mathcal{N}+1)^2$), \ref{lem-realvarLip} that
% ,  \ref{lem-normiequ} that 
\begin{equation*}
    \begin{aligned}
        \bigl\| e_\psi^{\frac{1}{2}} \bigr\| 
        &= \left\| g\left( \bigl| u^{\frac{1}{2}} \bigr|^2 \right)-g\left( \bigl| U^{\frac{1}{2}} \bigr|^2 \right) \right\| 
        \leq g_\sigma^1 \left\| \bigl|u^{\frac{1}{2}} \bigr|^2-\bigl| U^{\frac{1}{2}} \bigr|^2 \right\| \\
        &\leq g_\sigma^1 \left( \bigl\|u^{\frac{1}{2}}\bigr\|_{\infty} + \bigl\|U^{\frac{1}{2}}\bigr\|_{\infty}\right) \bigl\|e_u^{\frac{1}{2}}\bigr\| 
        \leq g_\sigma^1 \left(2\mathcal{N} + 1\right)\widetilde{C}\tau \left(\tau + h^2\right) \leq \widetilde{C}_0 \left(\tau^2 + h^2\right),
    \end{aligned}
\end{equation*}
and
\begin{equation*}
   \begin{aligned}
        \bigl\|e_\psi^{\frac{1}{2}} \bigr\|_{\infty} &= \left\|g\left( \bigl| u^{\frac{1}{2}} \bigr|^2 \right)-g\left( \bigl| U^{\frac{1}{2}} \bigr|^2 \right) \right\|_{\infty} 
        % \leq g_\sigma^1 \left\|\bigl| u^{\frac{1}{2}} \bigr|^2-\bigl| U^{\frac{1}{2}}\bigr|^2 \right\|_{\infty} \\
        \leq g_\sigma^1 \left(\bigl\|u^{\frac{1}{2}} \bigr\|_{\infty} + \bigl\| U^{\frac{1}{2}} \bigr\|_{\infty}\right)\bigl\|e_u^{\frac{1}{2}}\bigr\|_{\infty}
        \leq \widetilde{C}_0 \left(\tau + h^2\right).
   \end{aligned}
\end{equation*}
This completes the proof for \eqref{eq-miderrest}.
\end{proof}

Since the numerical solution $U^k$ in \eqref{eq-orierreq4}--\eqref{eq-orierreq6} lacks the mesh-independent bound in the maximum norm, we can not estimate the error terms $e_u^k$ for $1\leq k\leq K-1$.
To overcome this difficulty, we introduce suitable cut-off coefficients to formulate a modified error system analogous to the error system \eqref{eq-orierreq4}--\eqref{eq-orierreq6}.

Introduce the auxiliary error functions $\xi_u^k,\quad \xi_\psi^{k-1/2}$ for $1\leq k\leq K$, with $\xi_u^0 = e_u^0 = 0$, $\xi_\psi^{1/2} = e_\psi^{1/2}$. 
Let $\varepsilon$ be a fixed constant satisfying $0<\varepsilon<1$ and $ \mathcal{N}_{\star}>0 $ be a constant to be determined later. Define the cut-off coefficients by
\begin{equation}\label{eq-cutoffcoe}
    \lambda_u:=\min\Biggl\{\frac{1}{\max_{1\leq k\leq K}\|\xi_u^k\|_\infty+\varepsilon},1\Biggr\},\qquad 
    \lambda_\psi:=\min \Biggl\{\frac{\mathcal{N}_{\star}}{\max_{1\leq k\leq K} \bigl\| \xi_\psi^{k-\frac{1}{2}} \bigr\|_\infty+\varepsilon},1 \Biggr\}.
\end{equation}
Correspondingly, the modified numerical variables are denoted as
\begin{equation}\label{eq-modvarible}
    U_\lambda^k=u^k - \lambda_u\xi_u^k, \qquad \Psi_\lambda^{k-\frac{1}{2}} =\psi^{k-\frac{1}{2}}-\lambda_\psi\xi_\psi^{k-\frac{1}{2}}.
\end{equation}
With these preparations, we introduce the following modified error equations: 
\begin{equation}\label{eq-auxerreq1}
    \bi\delta_t \xi_u^{k+1}-(-\Delta_h)^{\frac{\alpha}{2}} \wh{\xi}_u^{\,k+\frac{1}{2}}=\mathrm{\Theta}_{u,1}^{k+\frac{1}{2}}+\mathrm{\Theta}_{u,2}^{k+\frac{1}{2}}+\mathrm{\Theta}_{u,3}^{k+\frac{1}{2}}, \qquad 0\leq k\leq K-1,
\end{equation}
% for $0\leq k\leq K-1$, 
where 
\begin{equation*}
    \mathrm{\Theta}_{u,1}^{k+\frac{1}{2}}=-\lambda_\psi \xi_\psi^{k+\frac{1}{2}}\wh{U}^{k+\frac{1}{2}}_\lambda, \qquad 
    \mathrm{\Theta}_{u,2}^{k+\frac{1}{2}}=-\psi^{k+\frac{1}{2}}\lambda_u\wh{\xi}_u^{\,k+\frac{1}{2}}, \qquad 
    \mathrm{\Theta}_{u,3}^{k+\frac{1}{2}}= \mathtt{r}_{u}^{k+\frac{1}{2}},   
\end{equation*}
and 
% for $1\leq k\leq K-1$, 
\begin{equation}\label{eq-auxerreq3}
    \xi_\psi^{k+\frac{1}{2}}+\xi_\psi^{k-\frac{1}{2}}
    =2\left(g\left(| u^k |^2\right)-g\left(| U_\lambda^k |^2\right)\right)+\mathtt{r}_\psi^k,\qquad 1\leq k\leq K-1.
\end{equation}

\begin{remark}
    For the fixed cut-off coefficients $(\lambda_u, \lambda_{\psi}) \in [0,1]^2$, the modified error system \eqref{eq-auxerreq1}--\eqref{eq-auxerreq3} is uniquely solvable, and its solution continuously depends on these coefficients. Therefore, the definition of cut-off coefficients induces a continuous self-mapping on $[0,1]^2$. Using Brouwer's fixed-point theorem, we can prove the existence of solutions to the modified error system.
\end{remark}

\begin{remark}
By the definitions of the cut-off coefficients in \eqref{eq-cutoffcoe}, the auxiliary error functions satisfy
\begin{equation}\label{eq-merrbd}
    \bigl\|\lambda_u\xi_u^{k}\bigr\|_\infty  \leq 1, \qquad
    \bigl\|\lambda_\psi\xi_\psi^{k-\frac{1}{2}} \bigr\|_\infty \leq \mathcal{N}_{\star}, \quad 1\leq k\leq K.
\end{equation}
According to the boundedness assumption \eqref{eq-exsolbound}, the modified variables \eqref{eq-modvarible} imply that
\begin{equation}\label{eq-mvarbd}
    \bigl\| U_\lambda^k \bigr\|_\infty  \leq \mathcal{N}+1, \qquad 
    \bigl\| \Psi_\lambda^{k-\frac{1}{2}} \bigr\|_\infty \leq \mathcal{N}+\mathcal{N}_{\star}, \quad 1\leq k\leq K.
\end{equation}
These bounds make it possible to estimate the auxiliary error functions in \eqref{eq-auxerreq1}--\eqref{eq-auxerreq3}. Once it is proven that both coefficients in \eqref{eq-cutoffcoe} are equal to one, the modified error system \eqref{eq-auxerreq1}--\eqref{eq-auxerreq3} coincides with the original error system \eqref{eq-orierreq4}--\eqref{eq-orierreq6}.

\end{remark}

\subsection{Closed estimates for modified errors}

Before deriving the modified error estimates via \eqref{eq-auxerreq1}--\eqref{eq-auxerreq3}, we set several specific positive constants: 
\begin{equation}\label{eq-specconst}
    \begin{aligned}
    &C_1:=\mathcal N+1, \quad  C_2=2(4\mathcal N+2)\left(g_\sigma^1+4\mathcal N^2g_\sigma^2\right),\quad  C_3:=C_1C_2+2\mathcal N, \quad C_4:=C_3+2C_1,\\
    &C_1^{\star}:=C_1C_2+3(\mathcal N+\mathcal N_{\star}), \qquad C_2^{\star}:= \max\left\{ 1+\mathcal N_{\star},\;C_1+C_3,\;C_2+C_1^\star\right\}.
    \end{aligned}
\end{equation}
These constants are independent of grid parameters, with the superscript $\star$ indicating dependence on the undetermined constant $\mathcal{N}_\star$. Furthermore, we also set some general positive constants $\{\widetilde{C}_q \}_{q=1}^{4}$ that are independent of grid parameters and $\mathcal{N}_\star$, whereas may be related to $\widetilde{C}$, $\mathcal{N}$, $T$, and the indicated derivatives $g_\sigma^1$, $g_\sigma^2$.
% Note that $g_\sigma^1$ and  $g_\sigma^2$
Note that $g_\sigma^1$ and  $g_\sigma^2$ are the constants introduced in Lemma \ref{lem-gfunLip}, its values depend on the parameter $\sigma$. In light of \eqref{eq-mvarbd} and \eqref{eq-exsolbound}, choosing $\sigma=(\mathcal{N}+1)^2$ ensures that the conditions of Lemma \ref{lem-gfunLip} are satisfied. For subsequent applications of Lemma \ref{lem-gfunLip}, we shall omit the choice of $\sigma$.

% whereas $\{\widetilde{C}^{\star}_q\}_{q=1}^{3}$ may be related to  $\mathcal{N}_\star$ but are still independent of the grid parameters.

% Finally, we shall repeatedly invoke Lemma \ref{lem-gfunLip} to handle $g(| u^k |^2)-g(| U_\lambda^k |^2)$ in \eqref{eq-auxerreq3}.
% In light of \eqref{eq-mvarbd} and \eqref{eq-exsolbound}, as long as $\sigma = (\mathcal{N}+1)^2$, then the condition of Lemma \ref{lem-gfunLip} is satisfied. For subsequent applications of Lemma \ref{lem-gfunLip}, the specific choice of $\sigma$ is omitted.

Now, we first provide an error estimate for $\|\xi_u^{n+1}\|$, which is controlled by $\|\xi_{\psi}^{k+1/2}\|$ for $0\leq k\leq n$. In addition, we also obtain the estimate for the first two time levels.

\begin{lemma}\label{lem-estimate1}
    Assume $u\in C^{4}\left([0,T],\mathcal{W}^{2+\alpha}(\mathbb R^2)\right)$, $g\in C^3(\mathbb{R};\mathbb{R})$. Then
    \begin{equation}
        % \|\xi_u^{n+1}\| \leq \tcb{(\mathcal{N}+1)}\tcr{\widetilde{C}_1}\tau  \sum_{k=0}^{n}\|\xi_{\psi}^{k+\frac{1}{2}}\| + \tcb{\widetilde{C}T}\tcr{\widetilde{C}_1} (\tau^2+h^2), \qquad 0\leq n\leq K-1,
        \bigl\| \xi_u^{n+1} \bigr\| 
        \leq C_1 \tau  \sum_{k=0}^{n} \bigl\|\xi_{\psi}^{k+\frac{1}{2}} \bigr\| + \widetilde{C}_1\left(\tau^2+h^2\right), \qquad 0\leq n\leq K-1,
    \end{equation}
    % where $\widetilde{C}_1$ is a positive constant that depends on $\widetilde{C}$, $\mathcal{N}$ and $T$.
    % As a special case, we get
    Moreover, at the first two time levels,
    \begin{equation}
        \bigl\| \xi_\psi^{m-\frac{1}{2}} \bigr\| + \bigl\|\delta_t \xi_u^{m} \bigr\| 
        \leq \widetilde{C}_1 \left( \tau^2+h^2 \right),\qquad 
        \bigl\|\xi_u^{m}\bigr\| 
        \leq \widetilde{C}_1 \tau\left( \tau^2+h^2 \right),\quad  m=1,2,
    \end{equation}
    where $\widetilde{C}_1$ is a positive constant that depends on $\widetilde{C}$, $\widetilde{C}_0$, $\mathcal{N}$, $T$ and $g_\sigma^1$.
\end{lemma}

\begin{proof}

Taking the inner product of \eqref{eq-auxerreq1} with $\wh\xi_u^{\,k+1/2}$ and then taking its imaginary part, we have
\begin{equation*}
\|\xi_u^{k+1}\|^2-\|\xi_u^k\|^2
=2\tau\impart\sum_{m=1}^3
 \langle\Theta_{u,m}^{k+1/2},\wh\xi_u^{\,k+1/2}\rangle.
\end{equation*}
Apply the Cauchy-Schwarz inequality, we get
    \begin{equation*}
        \bigl\| \xi_u^{k+1} \bigr\| - \bigl\| \xi_u^{k} \bigr\| \leq \tau \left( \bigl\| \mathrm{\Theta}_{u,1}^{k+\frac{1}{2}} \bigr\|+ \bigl\| \mathrm{\Theta}_{u,3}^{k+\frac{1}{2}} \bigr\|\right).
    \end{equation*}
    In light of \eqref{eq-mvarbd} and Lemma \ref{lem-truncerror}, we get 
    \begin{equation}\label{eq-fm1}
        \bigl\| \xi_u^{k+1} \bigr\| - \bigl\| \xi_u^{k} \bigr\| \leq  (\mathcal{N}+1) \tau \bigl\| \xi_{\psi}^{k+\frac{1}{2}} \bigr\| + \widetilde{C}\tau \left( \tau^2+h^2 \right),
    \end{equation}
    and summing this inequality from $k=0$ to $n$ yields 
    \begin{equation*}
        \bigl\| \xi_u^{n+1} \bigr\| 
        \leq (\mathcal{N}+1) \tau \sum_{k=0}^{n} \bigl\| \xi_{\psi}^{k+\frac{1}{2}} \bigr\| + \widetilde{C}T\left( \tau^2+h^2 \right) 
        \leq C_1 \tau \sum_{k=0}^{n} \bigl\| \xi_{\psi}^{k+\frac{1}{2}} \bigr\| + \widetilde{C}_1\left(\tau^2+h^2 \right).
    \end{equation*}

    In particular, for the equation \eqref{eq-fm1} with $k=0$, it follows from Lemma \ref{lem-halfrexest} that 
    \begin{equation}\label{eq-fm2}
        % \|\xi_u^{1}\|\leq \tau (\mathcal{N}+1) \|\xi_{\psi}^{\frac{1}{2}}\| + \widetilde{C}\tau(\tau^2+h^2) \leq \tcb{[(\mathcal{N}+1)\widetilde{C}_0+\widetilde{C}]}\tcr{\widetilde{C}_2}\tau(\tau^2+h^2).
        \bigl\| \xi_u^{1} \bigr\|\leq  (\mathcal{N}+1)\tau \bigl\| \xi_{\psi}^{\frac{1}{2}} \bigr\| + \widetilde{C}\tau\left(\tau^2+h^2\right) \leq \widetilde{C}_1\tau\left(\tau^2+h^2\right).
    \end{equation}
    Then, for the equation \eqref{eq-auxerreq3} with $k=1$, and using Lemmas \ref{lem-gfunLip}, \ref{lem-truncerror}, and \ref{lem-halfrexest}, together with \eqref{eq-mvarbd}, we have 
    \begin{equation*}
        \begin{aligned}
            \bigl\| \xi_\psi^{\frac{3}{2}} \bigr\| 
            &\leq \bigl\|\xi_\psi^{\frac{1}{2}}\bigr\| + \bigl\|\mathtt{r}_\psi^{1}\bigr\| + 2\norm{ g\left( |u^1|^2 \right)-g\left( |U_\lambda^1|^2 \right) }  \\
            &\leq \bigl\| \xi_\psi^{\frac{1}{2}} \bigr\| + \bigl\| \mathtt{r}_\psi^{1} \bigr\| + 2g_\sigma^1\left(\bigl\| u^1 \bigr\|_{\infty}+ \bigl\| U_\lambda^1 \bigr\|_{\infty}\right)\bigl\|\lambda_u \xi_u^1 \bigr\| \\
            % &\leq \tcb{\widetilde{C}_0}(\tau^2+h^2) + \tcb{\widetilde{C}}\tau^2 + \tcb{g_\sigma^1\left(4\mathcal{N}+ 2\right)[(\mathcal{N}+1)\widetilde{C}_0+\widetilde{C}]T}(\tau^2+h^2) \\
            % &\leq \tcb{\widetilde{C}_0}(\tau^2+h^2) + \tcb{\widetilde{C}}\tau^2 + g_\sigma^1\tcr{\left(4\mathcal{N}+ 2\right)\|\xi_u^1\|} 
            % \leq \tcr{\widetilde{C}_2}(\tau^2+h^2).
            &\leq \widetilde{C}_0 \left(\tau^2+h^2 \right) + \widetilde{C}\tau^2 + g_\sigma^1\left(4\mathcal{N}+ 2\right)\bigl\| \xi_u^1 \bigr\| 
            \leq \widetilde{C}_1 \left(\tau^2+h^2\right).
        \end{aligned}
    \end{equation*} 
    Similar to \eqref{eq-fm2}, we can derive from \eqref{eq-fm1} that
    \begin{equation*}
        \bigl\| \xi_u^{2} \bigr\| 
        \leq \bigl\| \xi_u^{1} \bigr\| +  (\mathcal{N}+1) \tau \bigl\| \xi_{\psi}^{\frac{3}{2}} \bigr\| + \widetilde{C} \tau \left( \tau^2+h^2 \right) \leq  \widetilde{C}_1\tau \left(\tau^2+h^2\right).
    \end{equation*}
    Finally, it is easy to obtain that
    \begin{equation*}
        \bigl\|\delta_t \xi_u^{1} \bigr\| = \tau^{-1} \bigl\|\xi_u^{1} \bigr\| \leq \widetilde{C}_1 \left(\tau^2+h^2 \right),\quad 
        \bigl\| \delta_t \xi_u^{2} \bigr\| = \tau^{-1} \bigl\|\xi_u^{2}-\xi_u^{1} \bigr\| \leq \tau^{-1} \left(\bigl\|\xi_u^{2} \bigr\| + \bigl\| \xi_u^{1} \bigr\| \right) \leq \widetilde{C}_1 \left(\tau^2+h^2\right).
    \end{equation*}
    Therefore, the proof is completed.
\end{proof}

% Next, we prove the estimate of $\bigl\|\xi_\psi^{n+\frac{1}{2}}\bigr\|$, which is related to $\norm{\xi_u^k}$ and $\norm{\delta_t\xi_u^{k}}$ for $2\leq k\leq n$.
Next, we estimate the error $\|\xi_\psi^{n+1/2}\|$, which is related to $\|\xi_u^k\|$ and $\|\delta_t\xi_u^{k}\|$ for $2\leq k\leq n$.

\begin{lemma}\label{lem-estimate2}
    % Under the same conditions as in Lemma \ref{lem-estimate1},
    Under the assumptions of Lemma \ref{lem-estimate1}, we obtain
    % there exists a positive constant $\widetilde{C}_3$ related to $\mathcal{N}$, $T$, $g_\sigma^1$, $g_\sigma^2$ and $\widetilde{C}_2$, such that
    \begin{equation}
        \bigl\|\xi_\psi^{n+\frac{1}{2}}\bigr\| 
        \leq C_2 \tau\sum_{k=2}^{n} \left(\norm{\xi_u^k} + \norm{\delta_t\xi_u^{k}}\right)+ \widetilde{C}_2 \left(\tau^2+h^2\right), \quad 2\leq n\leq K-1,
    \end{equation}
    where the positive constant $\widetilde{C}_2$ is related to $\widetilde{C}$, $T$ and $\widetilde{C}_1$.
\end{lemma}

\begin{proof}
    By subtracting \eqref{eq-auxerreq3} at $t = t_{k-1}$ from the equation at $t = t_k$, we have 
    \begin{equation}\label{eq-fm3}
        \xi_\psi^{k+\frac{1}{2}} - \xi_\psi^{k-\frac{3}{2}} = 2\mathfrak{G}^k  + \mathtt{r}_\psi^{k} -\mathtt{r}_\psi^{k-1}
    \end{equation}
for $2\leq k\leq K-1$, where
\begin{equation*}
    \mathfrak{G}^k = \mathrm{G}\bigl(|U_\lambda^{k-1}|^2,|U_\lambda^k|^2,|u^{k-1}|^2,|u^k|^2\bigr).
\end{equation*}
Set $\mathrm S^k:=\mathrm S(U_\lambda^{k-1},U_\lambda^k,u^{k-1},u^k)$, and applying Lemmas \ref{lem-gfunLip} and \ref{lem-realvarLip}, we derive
    \begin{equation*}
        \begin{aligned}
             \, \norm{\mathfrak{G}^k} 
            =& \, \norm{g\left(|U_\lambda^{k-1}|^2\right)-g\left(|U_\lambda^k|^2\right)- g\left(|u^{k-1}|^2\right)+g\left(|u^k|^2\right)} \\
            \leq& \, g_\sigma^{2} \bigl\| |u^{k-1}|^2-|u^k|^2 \bigr\|_{\infty} \bigl\||U_\lambda^k|^2-|u^{k}|^2  \bigr\| + \left(g_\sigma^{1}+g_\sigma^{2} \bigl\| |u^{k-1}|^2-|u^k|^2 \bigr\|_{\infty}\right) \norm{\mathrm{S}^{k}} \\
            \leq& \, 2\mathcal{N}g_\sigma^{2}\bigl\| u^{k-1}-u^k \bigr\|_{\infty} \cdot (2\mathcal{N} + 1)\bigl\| \xi_u^k \bigr\| + \left(g_\sigma^{1}+2\mathcal{N}g_\sigma^{2}\bigl\| u^{k-1}-u^k \bigr\|_{\infty}\right) \\
            &\cdot\left[2\bigl\| u^{k-1}-u^k \bigr\|_{\infty} \bigl\| \xi_u^k \bigr\| + \left(2\mathcal{N}+2+\bigl\| u^{k-1}-u^k \bigr\|_{\infty}\right) \norm{\xi_u^{k}-\xi_u^{k-1}}\right]  \\
        \end{aligned}
    \end{equation*}
    According to \eqref{eq-exsolbound}, it can be observed that 
    \begin{equation*}
        \|u^{k-1}-u^k\|_{\infty} \leq \min\left\{\tau\|\partial_t u\|_{L^{\infty}}, \|u^{k-1}\|_{\infty}+\|u^k\|_{\infty}\right\} \leq \min\left\{\mathcal{N}\tau,2\mathcal{N}\right\}.
    \end{equation*}
    Thus, 
    \begin{equation*}
        \begin{aligned}
            \norm{\mathfrak{G}^k} 
            \leq& \left[(4\mathcal{N}+2)\mathcal{N}^2 g_{\sigma}^2\right]\tau\norm{\xi_u^k} + \left[\left(g_{\sigma}^1 + 4\mathcal{N}^2 g_{\sigma}^2\right)\cdot 2\mathcal{N}\right]\tau\norm{\xi_u^k} \\
            &+ \left[\left(g_{\sigma}^1 + 4\mathcal{N}^2 g_{\sigma}^2\right)\cdot \left(4\mathcal{N}+2\right)\right]\tau\norm{\delta_t\xi_u^{k}} \\
            \leq & (4\mathcal{N}+2)\bigl(g_{\sigma}^1+4\mathcal{N}^2 g_{\sigma}^2\bigr) \tau \left(\norm{\xi_u^k} + \norm{\delta_t\xi_u^{k}}\right) \\
            =& \frac{C_2}{2} \tau\left(\norm{\xi_u^k} + \norm{\delta_t\xi_u^{k}}\right) .
        \end{aligned}
    \end{equation*}
    % where $c_1$ and $c_2$ are constants related to $g_{\sigma}^1$, $g_{\sigma}^2$, $\mathcal{N}$ and $T$.
    Taking the inner product of \eqref{eq-fm3} with $\xi_\psi^{k+1/2} + \xi_\psi^{k-3/2}$ and applying Lemma \ref{lem-truncerror}, we obtain
    % Now, by taking an inner product of \eqref{eq-fm3} with $\xi_\psi^{k+\frac{1}{2}} + \xi_\psi^{k-\frac{3}{2}}$, and applying Lemma \ref{lem-truncerror}, we obtain
    \begin{equation}\label{eq-fm7}
        \begin{aligned}
            \bigl\| \xi_\psi^{k+\frac{1}{2}} \bigr\| - \bigl\|\xi_\psi^{k-\frac{3}{2}} \bigr\| 
            \leq& 2\norm{\mathfrak{G}^k} + \bigl\| \mathtt{r}_\psi^{k} -\mathtt{r}_\psi^{k-1} \bigr\|  \\
            % \leq& \tcb{2\max\{g_\sigma^{2}\mathcal{N},g_\sigma^{1}+g_\sigma^{2}\mathcal{N}T\}}\left(\norm{\xi_u^k} + \norm{\delta_t\xi_u^{k}}\right) + \widetilde{C}\tau^3 \\
            % \leq& \tcr{\widetilde{C}_3}\tau \left(\norm{\xi_u^k} + \norm{\delta_t\xi_u^{k}}\right) + \widetilde{C}\tau^3,
            \leq& C_2 \tau \left(\norm{\xi_u^k} + \norm{\delta_t\xi_u^{k}}\right) + \widetilde{C}\tau^3,
        \end{aligned}
    \end{equation}
    and after summing over $k$ from $2$ to $n$ ($2\leq n\leq K-1$), we get
    \begin{equation*}
        \begin{aligned}
            \bigl\| \xi_\psi^{n+\frac{1}{2}} \bigr\| + \bigl\|\xi_\psi^{n-\frac{1}{2}} \bigr\|  
             \leq& \bigl\| \xi_\psi^{\frac{1}{2}} \bigr\| + \bigl\| \xi_\psi^{\frac{3}{2}} \bigr\|+ C_2 \tau\sum_{k=2}^{n} \left(\norm{\xi_u^k} + \norm{\delta_t\xi_u^{k}}\right) + \widetilde{C}T\tau^2  \\
            % \leq& \tcb{2C}\tcr{\widetilde{C}_3}\tau\sum_{k=2}^{n} \left(\norm{\xi_u^k} + \norm{\delta_t\xi_u^{k}}\right)+\tcb{\max\{C,\widetilde{C}T\}} \tcr{\max\{\widetilde{C}_2,\widetilde{C}T\}}\tcr{\widetilde{C}_4}(\tau^2+h^2), \\
            \leq& C_2 \tau\sum_{k=2}^{n} \left(\norm{\xi_u^k} + \norm{\delta_t\xi_u^{k}}\right) + \widetilde{C}_2 \left(\tau^2+h^2\right),
        \end{aligned}
    \end{equation*}
    where the Lemma \ref{lem-estimate1} is used. This finishes the proof.
\end{proof}

Below, we build an estimate for $\bigl\|\delta_t \xi_u^{n+1}\bigr\|$, which is crucial for deriving the final result.

\begin{lemma}\label{lem-estimate3}
    Under the assumptions of Lemma \ref{lem-estimate1}, and for $\tau(1+\mathcal N_\star)\leq1$, we have
    % we can deduce that
    \begin{equation}
        \bigl\| \delta_t \xi_u^{n+1} \bigr\| 
        \leq C_1^{\star} \tau \sum_{k=2}^{n}\norm{\delta_t\xi_u^{k+1}} +  C_3\tau \sum_{k=2}^{n} \bigl( \bigl\| \xi_\psi^{k-\frac{3}{2}} \bigr\| + \norm{\xi_u^{k+1}} \bigr)+ \widetilde{C}_3\left(\tau^2+h^2\right), \quad 2\leq n\leq K-1,
    \end{equation}
    where $\widetilde{C}_3$ is a positive constant that depends on $\mathcal{N}$, $\widetilde{C}_1$, $T$, $\widetilde{C}$ and $C_1$.
\end{lemma}

\begin{proof}
    By subtracting \eqref{eq-auxerreq1} at $t = t_{n-3/2}$ from the equation at $t = t_{n+1/2}$, we get
    \begin{equation}\label{eq-fm4}
        \delta_t \xi_u^{n+1}-\delta_t \xi_u^{n-1}+\frac{\bi}{2}(-\Delta_h)^{\frac{\alpha}{2}} \left(\xi_u^{n+1}+\xi_u^{n}-\xi_u^{n-1}-\xi_u^{n-2} \right) = \sum_{m=1}^{3} \left( \bi\mathrm{\Theta}_{u,m}^{n-\frac{3}{2}} - \bi\mathrm{\Theta}_{u,m}^{n+\frac{1}{2}} \right)
    \end{equation}
    for $2\leq n\leq K-1$, and rewrite 
    \begin{equation*}
        \xi_u^{n+1}+\xi_u^{n}-\xi_u^{n-1}-\xi_u^{n-2} = \tau\left(\delta_t\xi_u^{n+1}+2\delta_t\xi_u^{n}+\delta_t\xi_u^{n-1} \right).
    \end{equation*}
    Recalling the operators 
    $\mathcal{P}_h = \mathcal{I}_h+\frac{\bi \tau}{2}(-\Delta_h)^{\frac{\alpha}{2}}$, 
    $\mathcal{Q}_h = \mathcal{I}_h-\frac{\bi \tau}{2}(-\Delta_h)^{\frac{\alpha}{2}}$ and $\mathcal{B}_h = \mathcal{P}_h^{-1} \mathcal{Q}_h$, thereby \eqref{eq-fm4} can be written as
    \begin{equation*}
        \delta_t \xi_u^{n+1} = \left(\mathcal{B}_h - \mathcal{I}_h\right)\delta_t \xi_u^{n} + \mathcal{B}_h \delta_t \xi_u^{n-1} + \mathcal{P}_h^{-1}\sum_{m=1}^{3}\mathrm{\Upsilon}_{m}^{n},
    \end{equation*}
    where $\mathrm{\Upsilon}_{m}^{n} := \bi\mathrm{\Theta}_{u,m}^{n-3/2} - \bi\mathrm{\Theta}_{u,m}^{n+1/2}$. Thus, it follows from Lemma \ref{lem-opeequality} that
    \begin{equation}\label{eq-fm5}
        \begin{aligned}
            \delta_t \xi_u^{n+1} =& \frac{1}{2}\left[(-1)^{n+1}\mathcal{I}_h + \mathcal{B}_h^{n}\right]\mathcal{P}_h \delta_t \xi_u^{2} + \frac{1}{2}\left[(-1)^{n}\mathcal{B}_h + \mathcal{B}_h^{n}\right]\mathcal{P}_h \delta_t \xi_u^{1} \\ 
            &+ \frac{1}{2}\sum_{k=2}^{n}\left[(-1)^{n-k}\mathcal{I}_h + \mathcal{B}_h^{n-k+1} \right]\sum_{m=1}^{3}\mathrm{\Upsilon}_{m}^{k},
        \end{aligned}
    \end{equation} 
    and using Lemma \ref{lem-openorm}, we can deduce from \eqref{eq-fm5} that 
    \begin{equation}\label{eq-fm6}
        \bigl\|\delta_t \xi_u^{n+1} \bigr\| \leq  \bigl\|\mathcal{P}_h \delta_t \xi_u^{2} \bigr\| + \bigl\|\mathcal{P}_h \delta_t \xi_u^{1} \bigr\| + \sum_{k=2}^{n}\sum_{m=1}^{3} \bigl\| \mathrm{\Upsilon}_{m}^{k} \bigr\|.
    \end{equation}

    Now, we estimate the terms on the right-hand side of \eqref{eq-fm6} sequentially. For the equation \eqref{eq-auxerreq1} with $k=0,1$, it can be derived that
    \begin{equation*}
        \mathcal{P}_h\delta_t \xi_u^{1} = -\sum_{m=1}^{3}\bi\mathrm{\Theta}_{u,m}^{\frac{1}{2}},\qquad
        \mathcal{P}_h\delta_t \xi_u^{2} = - 2\left(\mathcal{P}_h -\mathcal{I}_h \right)\delta_t\xi_u^{1}-\sum_{m=1}^{3}\bi\mathrm{\Theta}_{u,m}^{\frac{3}{2}}.
    \end{equation*}
    Combining \eqref{eq-exsolbound} and \eqref{eq-mvarbd} with Lemmas \ref{lem-truncerror} and \ref{lem-estimate1}, we obtain
    % it can be seen that 
    \begin{equation*}
        % \tcb{
        % \| \mathrm{\Theta}_{u,1}^{q-\frac{1}{2}} \|\leq (\mathcal{N}+1)C \tcr{C_1\widetilde{C}_2}\left(\tau^2+h^2\right) , \quad \|\mathrm{\Theta}_{u,2}^{q-\frac{1}{2}}\|\leq \mathcal{N}CT \tcr{\mathcal{N}\widetilde{C}_2T}\left(\tau^2+h^2\right), \quad \|\mathrm{\Theta}_{u,3}^{q-\frac{1}{2}}\|\leq \widetilde{C}\left(\tau^2+h^2\right), \quad q=1,2 }
        \bigl\| \mathrm{\Theta}_{u,1}^{p-\frac{1}{2}} \bigr\|\leq (\mathcal{N}+1)\widetilde{C}_1\left(\tau^2+h^2\right), \quad 
        \bigl\| \mathrm{\Theta}_{u,2}^{p-\frac{1}{2}} \bigr\|\leq \mathcal{N}\widetilde{C}_1 T\left(\tau^2+h^2\right), \quad 
        \bigl\| \mathrm{\Theta}_{u,3}^{p-\frac{1}{2}} \bigr\|\leq \widetilde{C}\left(\tau^2+h^2\right), \quad p=1,2.
    \end{equation*}
    % \begin{equation}
    %     \sum_{m=1}^{3}\|\mathrm{\Theta}_{u,m}^{q-\frac{1}{2}}\| \leq \tcb{C}\tcr{\widetilde{C}_4}\left(\tau^2+h^2\right), \qquad   q = 1,2.
    % \end{equation}
    Therefore,
    \begin{equation}\label{eq-fm11}
        \begin{aligned}
            &\bigl\| \mathcal{P}_h \delta_t \xi_u^{1} \bigr\| \leq \sum_{m=1}^{3} \bigl\| \mathrm{\Theta}_{u,m}^{\frac{1}{2}} \bigr\| \leq \widetilde{C}_3\left(\tau^2+h^2\right), \\
            &\bigl\| \mathcal{P}_h \delta_t \xi_u^{2} \bigr\| \leq 2\left( \bigl\| \mathcal{P}_h \delta_t \xi_u^{1} \bigr\|+ \bigl\| \delta_t \xi_u^{1} \bigr\|\right) + \sum_{m=1}^{3}\bigl\|\mathrm{\Theta}_{u,m}^{\frac{3}{2}} \bigr\| \leq \widetilde{C}_3\left(\tau^2+h^2\right).
        \end{aligned}
    \end{equation}

    Next, we analyze $\|\mathrm{\Upsilon}_m^k\|$ for $m=1,2,3$.
    Regarding $\|\mathrm{\Upsilon}_1^k\|$, we split it into  two parts,
    \begin{equation*}
        \begin{aligned}
            \bigl\| \mathrm{\Upsilon}_1^k \bigr\| 
            &= \bigl\| \lambda_\psi \xi_\psi^{k+\frac{1}{2}}\wh{U}^{k+\frac{1}{2}}_\lambda -\lambda_\psi \xi_\psi^{k-\frac{3}{2}}\wh{U}^{k-\frac{3}{2}}_\lambda \bigr\| \\
            &\leq \bigl\| \bigl(\lambda_\psi \xi_\psi^{k+\frac{1}{2}} -\lambda_\psi \xi_\psi^{k-\frac{3}{2}}\bigr)\wh{U}^{k+\frac{1}{2}}_\lambda \bigr\| + \bigl\|\lambda_\psi \xi_\psi^{k-\frac{3}{2}}\bigl(\wh{U}^{k+\frac{1}{2}}_\lambda-\wh{U}^{k-\frac{3}{2}}_\lambda\bigr) \bigr\| \\
            &:= \mathrm{\Pi}_1 + \mathrm{\Pi}_2 .
        \end{aligned}
    \end{equation*}
    Using \eqref{eq-exsolbound}, \eqref{eq-merrbd}--\eqref{eq-mvarbd}, \eqref{eq-fm3} and \eqref{eq-fm7}, we obtain
    \begin{equation}\label{eq-fm16}
        \begin{aligned}
            \mathrm{\Pi}_1
            &\leq \bigl\| \wh{U}^{k+\frac{1}{2}}_\lambda \bigr\|_{\infty} \bigl\| \xi_\psi^{k+\frac{1}{2}} - \xi_\psi^{k-\frac{3}{2}}\bigr\| \\
            &\leq \bigl\| \wh{U}^{k+\frac{1}{2}}_\lambda\bigr\|_{\infty} \left(2\norm{\mathfrak G^k} + \bigl\| \mathtt{r}_\psi^{k} -\mathtt{r}_\psi^{k-1} \bigr\| \right) \\
            % &\leq \tcb{2C(\mathcal{N}+1)} \tcr{C_1\widetilde{C}_3}\tau \left(\norm{\xi_u^k} + \norm{\delta_t\xi_u^{k}}\right) + \tcb{\widetilde{C}(\mathcal{N}+1)}\tcr{C_1\widetilde{C}}\tau(\tau^2 + h^2),
            % &\leq 2(\mathcal{N}+1) \tcr{C_1\widetilde{C}_3}\tau \left(\norm{\xi_u^k} + \norm{\delta_t\xi_u^{k}}\right) + \tcb{\widetilde{C}(\mathcal{N}+1)}\tcr{C_1\widetilde{C}}\tau(\tau^2 + h^2),
            &\leq C_1 \left[C_2\tau \left(\norm{\xi_u^k} + \norm{\delta_t\xi_u^{k}}\right) + \widetilde{C}\tau^3\right], 
        \end{aligned}
    \end{equation}
    and
    \begin{equation}\label{eq-fm17}
        \begin{aligned}
            \mathrm{\Pi}_2
            &= \left\| \lambda_\psi \xi_\psi^{k-\frac{3}{2}}\Bigl[\wh{u}^{k+\frac{1}{2}}-\wh{u}^{k-\frac{3}{2}} - \bigl(\lambda_u\wh{\xi}_u^{\,k+\frac{1}{2}} - \lambda_u\wh{\xi}_u^{\,k-\frac{3}{2}} \bigr)\Bigr] \right\|  \\
            &\leq  \bigl\|\wh{u}^{k+\frac{1}{2}}-\wh{u}^{k-\frac{3}{2}} \bigr\|_{\infty} \bigl\|\xi_\psi^{k-\frac{3}{2}}\bigr\| + \bigl\|\lambda_\psi \xi_\psi^{k-\frac{3}{2}}\bigr\|_{\infty} \bigl\|\lambda_u\wh{\xi}_u^{\,k+\frac{1}{2}} - \lambda_u\wh{\xi}_u^{\,k-\frac{3}{2}}\bigr\| \\
            % &\leq 2\tau\norm{\partial_t u}_{L^\infty}\norm{\xi_\psi^{k-\frac{3}{2}}} + \frac{\mathcal{N}_{\star}\tau}{2} \norm{\delta_t\xi_u^{\,k+1} + 2\delta_t\xi_u^{\,k}+\delta_t\xi_u^{\,k-1}} \\
            % &\leq \tcb{2\mathcal{N}}\tau \norm{\xi_\psi^{k-\frac{3}{2}}} + \tcb{\mathcal{N}_{\star}\frac{\tau}{2}} \left(\norm{\delta_t\xi_u^{\,k+1}} + 2\norm{\delta_t\xi_u^{\,k}} +\norm{\delta_t\xi_u^{\,k-1}}\right).\\
            &\leq 2\mathcal{N}\tau \bigl\| \xi_\psi^{k-\frac{3}{2}} \bigr\| + \mathcal{N}_{\star}\tau\left(\norm{\delta_t\xi_u^{k+1}} + \norm{\delta_t\xi_u^{k}} +\norm{\delta_t\xi_u^{k-1}} \right).
        \end{aligned}
    \end{equation}
    So we further get
    \begin{equation}\label{eq-fm8}
        \begin{aligned}
            \bigl\| \mathrm{\Upsilon}_1^k \bigr\| 
            \leq& C_1C_2 \tau \left(\bigl\|\xi_u^k \bigr\|+  \norm{\delta_t\xi_u^{k}}\right) + \mathcal{N}_{\star}\tau \left(\norm{\delta_t\xi_u^{k+1}} +\norm{\delta_t\xi_u^{k}}+ \norm{\delta_t\xi_u^{k-1}} \right)  \\
            & + 2\mathcal{N}\tau \bigl\| \xi_\psi^{k-\frac{3}{2}}\bigr\| + C_1\widetilde{C}\tau\left(\tau^2 + h^2\right).
            %  + \widetilde{C}_1^{\star}\tau\left(\norm{\delta_t\xi_u^{\,k+1}} + \norm{\delta_t\xi_u^{\,k}} +\norm{\delta_t\xi_u^{\,k-1}} \right) +  \widetilde{C}_3\tau \left(\tau^2 + h^2\right). \\
        \end{aligned}
    \end{equation}
    As for $\bigl\|\mathrm{\Upsilon}_2^k \bigr\|$, it can be also seen from \eqref{eq-exsolbound} that
    \begin{equation}\label{eq-fm9}
        \begin{aligned}
            \bigl\| \mathrm{\Upsilon}_2^k \bigr\| 
            &= \bigl\| \psi^{k+\frac{1}{2}}\lambda_u\wh{\xi}_u^{\,k+\frac{1}{2}} -\psi^{k-\frac{3}{2}}\lambda_u\wh{\xi}_u^{\,k-\frac{3}{2}} \bigr\| \\
            % &\leq  \norm{\left(\psi^{k+\frac{1}{2}} - \psi^{k-\frac{3}{2}}\right)\lambda_u\wh{\xi}_u^{\,k+\frac{1}{2}}} + \norm{\psi^{k-\frac{3}{2}}\left(\lambda_u\wh{\xi}_u^{\,k+\frac{1}{2}} - \lambda_u\wh{\xi}_u^{\,k-\frac{3}{2}}\right)} \\
            &\leq  \bigl\| \psi^{k+\frac{1}{2}} - \psi^{k-\frac{3}{2}}\bigr\|_{\infty} \bigl\| \lambda_u\wh{\xi}_u^{\,k+\frac{1}{2}}\bigr\| + \bigl\| \psi^{k-\frac{3}{2}}\bigr\|_{\infty} \bigl\| \lambda_u\wh{\xi}_u^{\,k+\frac{1}{2}} - \lambda_u\wh{\xi}_u^{\,k-\frac{3}{2}}\bigr\| \\
            &\leq  \tau \norm{\partial_t \psi}_{L^\infty} \norm{\xi_u^{k+1}+\xi_u^{k}} + \frac{\tau}{2} \bigl\|\psi^{k-\frac{3}{2}}\bigr\|_{\infty} \norm{\delta_t \xi_u^{k+1} + 2\delta_t \xi_u^{k} + \delta_t \xi_u^{k-1}} \\
            &\leq \mathcal{N}\tau \left(\norm{\xi_u^{k+1}}+\norm{\xi_u^{k}} + \norm{\delta_t \xi_u^{k+1}} + \norm{\delta_t \xi_u^{k}} + \norm{\delta_t \xi_u^{k-1}}\right).
        \end{aligned}
    \end{equation}
    With respect to $\|\mathrm{\Upsilon}_3^k\|$, it can be known from Lemma \ref{lem-truncerror} that 
    \begin{equation}\label{eq-fm10}
        \bigl\| \mathrm{\Upsilon}_3^k \bigr\| 
        \leq \bigl\| \mathtt{r}_{u}^{k+\frac{1}{2}}-\mathtt{r}_{u}^{k-\frac{3}{2}}\bigr\| 
        \leq \bigl\| \mathtt{r}_{u}^{k+\frac{1}{2}}-\mathtt{r}_{u}^{k-\frac{1}{2}}\bigr\| + \bigl\| \mathtt{r}_{u}^{k-\frac{1}{2}}-\mathtt{r}_{u}^{k-\frac{3}{2}} \bigr\|
        \leq \widetilde{C}\tau\left(\tau^2+h^2\right).
    \end{equation}
    Then, Combining \eqref{eq-fm8}--\eqref{eq-fm10}, we obtain
    \begin{equation}\label{eq-threesum1}
        \begin{aligned}
        \sum_{m=1}^{3} \bigl\| \mathrm{\Upsilon}_{m}^{k} \bigr\|
        \leq &  \tau\left[ \mathcal{N}\norm{\xi_u^{k+1}} +(C_1C_2+\mathcal{N})\norm{\xi_u^{k}}+2\mathcal{N}\bigl\| \xi_\psi^{k-\frac{3}{2}}\bigr\|\right]+C_1C_2\tau \norm{\delta_t\xi_u^{k}} \\
        &+(\mathcal{N}_{\star}+\mathcal{N})\tau(\norm{\delta_t\xi_u^{k+1}}+\norm{\delta_t\xi_u^{k}}+\norm{\delta_t\xi_u^{k-1}})+\left(C_1\widetilde{C}+\widetilde{C}\right)T\left(\tau^2+h^2\right).
        \end{aligned}
    \end{equation}

    Finally, substituting \eqref{eq-fm11} and \eqref{eq-threesum1} into \eqref{eq-fm6} and applying Lemma \ref{lem-estimate1}, we have 
    \begin{equation*}
        \begin{aligned}
        \bigl\| \delta_t \xi_u^{n+1} \bigr\| 
        \leq & (C_1C_2+2\mathcal{N})\tau\sum_{k=2}^{n}\norm{\xi_u^{k+1}} + 2\mathcal{N}\tau\sum_{k=2}^{n}\bigl\| \xi_\psi^{k-\frac{3}{2}}\bigr\| + (C_1C_2+\mathcal{N})\norm{\xi_u^{2}} \\
        & + (C_1C_2+3(\mathcal{N}+\mathcal{N}_{\star}))\tau \sum_{k=2}^{n} \norm{\delta_t\xi_u^{k+1}}  + (\mathcal{N}+\mathcal{N}_{\star}) \tau \norm{\delta_t\xi_u^{1}} \\
        &+ (C_1C_2+2(\mathcal{N}+\mathcal{N}_{\star})) \tau \norm{\delta_t\xi_u^{2}} + \widetilde{C}_3\left(\tau^2+h^2\right)\\
        \leq & C_3\tau\sum_{k=2}^{n}\left(\norm{\xi_u^{k+1}} + \bigl\| \xi_\psi^{k-\frac{3}{2}}\bigr\|\right) + C_1^{\star}\tau \sum_{k=2}^{n} \norm{\delta_t\xi_u^{k+1}} \\ 
        &+ \mathcal{N}_{\star} \tau \left(\norm{\delta_t\xi_u^{1}} + 2\norm{\delta_t\xi_u^{2}}\right) + \widetilde{C}_3\left(\tau^2+h^2\right).
        \end{aligned}
    \end{equation*}
    We further apply the condition $\tau(1+\mathcal{N}_{\star})\leq 1$ and Lemma \ref{lem-estimate1} to the above inequality, and obtain
    \begin{equation}
        \bigl\| \delta_t \xi_u^{n+1} \bigr\| \leq  C_3\tau\sum_{k=2}^{n}\left(\norm{\xi_u^{k+1}} + \bigl\| \xi_\psi^{k-\frac{3}{2}}\bigr\|\right) + C_1^{\star}\tau \sum_{k=2}^{n} \norm{\delta_t\xi_u^{k+1}} +\widetilde{C}_3\left(\tau^2+h^2\right),
    \end{equation}
    which means that the proof has been completed.
\end{proof}

The above lemmas enable us to derive the following closed estimate for the modified system.

\begin{theorem}\label{thm-1}
    Assume that $u\in C^{4}\left([0,T],\mathcal{W}^{2+\alpha}(\mathbb R^2)\right)$ and $g\in C^3(\mathbb{R};\mathbb{R})$. Then, when $\tau \leq \tau_0 :=1/(2C_2^\star)$, we can show that
    \begin{equation}\label{eq-fm15}
       \max_{1\leq k\leq K}\Bigl( \bigl\| \xi_u^{k} \bigr\| + \bigl\|\xi_\psi^{k-\frac{1}{2}} \bigr\| + \bigl\|\delta_t \xi_u^{k} \bigr\|\Bigr) \leq C_3^{\star} \left(\tau^2+h^2\right),
    \end{equation}
    where $C_3^{\star}$ is a positive constant related to $\mathcal{N}_\star$, but is independent of $\tau$ and $h$.
\end{theorem}

\begin{proof}
    Let $\chi^{n}:= \bigl\| \xi_u^{n} \bigr\| + \bigl\| \xi_\psi^{n-1/2} \bigr\| + \bigl\| \delta_t \xi_u^{n} \bigr\|$ for $1\leq  n\leq K$. When $n=1,2$, it can be directly obtained from Lemma \ref{lem-estimate1} that 
    \begin{equation}\label{eq-fm13}
        \chi^{1} + \chi^{2}\leq C \left(\tau^2+h^2\right), 
    \end{equation}
Then, according to Lemmas \ref{lem-estimate1}--\ref{lem-estimate3}, we have 
\begin{equation}\label{eq-fm12}
    \begin{aligned}
        \chi^{n+1}
        \leq & \max\left\{C_3,C_1,C_1^{\star}\right\}\tau\chi^{n+1} + \max\left\{C_2+C_3, C_1+C_3, C_2+C_1^\star\right\}\tau\sum_{k=2}^n \chi^{k} \\
        &+ \widetilde{C}_4(\tau^2+h^2)\\
        \leq & C_2^{\star}\tau\chi^{n+1} + C_2^{\star}\tau\sum_{k=2}^n\chi^k + \widetilde{C}_4(\tau^2+h^2),\qquad 2\leq n\leq K-1.
    \end{aligned}
\end{equation}
Recalling the definition of $C_2^\star$ in \eqref{eq-specconst}, $C_2^\star\geq 1+\mathcal{N}_\star$ is hold. If we take $\tau \leq \tau_0 :=1/(2C_2^\star)$, the hypothesis of Lemma \ref{lem-estimate3} holds and the first term on the right-hand side of \eqref{eq-fm12} can be absorbed into the left-hand side.
With the aid of the discrete Gr\"{o}nwall's inequality, we can conclude from \eqref{eq-fm12} that
    \begin{equation}\label{eq-fm14}
        \chi^{n+1} \leq C_3^{\star} \left(\tau^2+h^2\right), \qquad 2\leq n\leq K-1,
    \end{equation}
    Finally, combining \eqref{eq-fm13} and \eqref{eq-fm14} gives \eqref{eq-fm15}, which means the proof has been completed.   
\end{proof}

\subsection{Uniform boundedness and optimal convergence}

Building on the modified error estimates established in the previous subsection, we now prove the unconditional optimal-order convergence of the RCND scheme.

Below, we first establish the uniform bounds of the modified errors $\xi_u^{k}$ and $\xi_\psi^{k-1/2}$
and determine the value of the undetermined constant $\mathcal{N}_\star$.

\begin{theorem}\label{thm-2}
    Assume that the conditions of Theorem \ref{thm-1} hold, and let $\xi_u^{k}\in \mathsf{C}_h^0$ and  $\xi_\psi^{k-1/2} \in \mathsf{R}_h^0$ be the solutions of the modified error equations \eqref{eq-auxerreq1}--\eqref{eq-auxerreq3}.
    Then, for sufficiently small $\tau$ and $h$, 
    there exists a constant $C^{\dagger}_{\alpha}>0$, depending on $\mathcal{N}$, $T$, $\alpha$, $g_{\sigma}^1$ and $g_{\sigma}^2$, but independent of $\tau$, $h$, and $\mathcal N_\star$, such that 
    \begin{equation}\label{eq-uniformbd}
        \max_{1\leq k\leq K} \norm{\xi_u^{k}}_{\infty} \leq 1-\varepsilon, \qquad
        \max_{1\leq k\leq K} \bigl\| \xi_\psi^{k-\frac{1}{2}}\bigr\|_{\infty} \leq \mathcal{N}_\star-\varepsilon, 
    \end{equation}
    where $\varepsilon$ is a fixed small constant defined in \eqref{eq-cutoffcoe}, and $\mathcal{N}_\star = 2C^{\dagger}_{\alpha}+\varepsilon$. 
\end{theorem}

\begin{proof}
    We split the proof into two complementary cases. In the case of $\tau > h$, we first apply \eqref{eq-auxerreq1}, Lemma \ref{lem-truncerror} and Theorem \ref{thm-1} to derive
    \begin{equation}\label{eq-fm18}
        \begin{aligned}
            \bigl\| (-\Delta_h)^{\frac{\alpha}{2}}\wh{\xi}_u^{\,k+\frac{1}{2}} \bigr\| 
            =& \Bigl\|\bi\delta_t\xi_u^{k+1}-\sum_{m=1}^3\Theta_{u,m}^{k+\frac12}\Bigr\|  \\
            \leq& C \left(\norm{\delta_t \xi_u^{k+1}} + \bigl\| \xi_\psi^{k+\frac{1}{2}} \bigr\| + \bigl\| \wh{\xi}_u^{\,k+\frac{1}{2}}\bigr\| + \bigl\|\mathtt{r}_u^{k+\frac{1}{2}}\bigr\| \right) \\
            \leq& C^\star_3 \left(\tau^2+h^2\right). 
        \end{aligned}
    \end{equation}
    Then, it follows from Lemma \ref{lem-uniform_maxest} that 
    \begin{equation}\label{eq-fm19}
        \bigl\| \wh{\xi}_u^{\,k+\frac{1}{2}}\bigr\|_{\infty} 
        \leq  C_{\alpha} \left(\bigl\|\wh{\xi}_u^{\,k+\frac{1}{2}}\bigr\| + \bigl\| (-\Delta_h)^{\frac{\alpha}{2}}\wh{\xi}_u^{\,k+\frac{1}{2}} \bigr\| \right) 
        \leq \widetilde{C}^\star_\alpha \left(\tau^2+h^2\right),
    \end{equation}
    where the constant $\widetilde{C}^\star_\alpha$ depends not only on $C^\star_3$, but also on a constant related to $\alpha$. Return to the estimate of $\mathrm{\Pi}_2$ in \eqref{eq-fm17}, we have
    \begin{equation*}
        \begin{aligned}
            \mathrm{\Pi}_2
            &\leq  \bigl\| \wh{u}^{k+\frac{1}{2}}-\wh{u}^{k-\frac{3}{2}}\bigr\|_{\infty} \bigl\|\xi_\psi^{k-\frac{3}{2}}\bigr\| +  \bigl\| \wh{\xi}_u^{\,k  +\frac{1}{2}} - \wh{\xi}_u^{\,k-\frac{3}{2}}\bigr\|_{\infty} \bigl\|\lambda_\psi \xi_\psi^{k-\frac{3}{2}}\bigr\| \\
            &\leq 2\mathcal{N}\tau \bigl\|\xi_\psi^{k-\frac{3}{2}}\bigr\| + 2\widetilde{C}^\star_\alpha (\tau^2+h^2) \bigl\| \xi_\psi^{k-\frac{3}{2}} \bigr\| 
            \leq 4\mathcal{N}\tau \bigl\|\xi_\psi^{k-\frac{3}{2}}\bigr\|,
        \end{aligned}
    \end{equation*}
    where the last inequality holds for sufficiently small $\tau$ and $h$. Then, by combining the above estimate for $\mathrm{\Pi}_2$ with the estimate of $\mathrm{\Pi}_1$ in \eqref{eq-fm16}, the estimate $\mathrm{\Upsilon}_1^k $ can be updated as 
    \begin{equation*}
        \bigl\| \mathrm{\Upsilon}_1^k \bigr\| 
        \leq C_1C_2 \tau \left(\bigl\|\xi_u^k \bigr\|+  \norm{\delta_t\xi_u^{k}}\right)
         + 4\mathcal{N}\tau \bigl\| \xi_\psi^{k-\frac{3}{2}}\bigr\| + C_1\widetilde{C}\tau\left(\tau^2 + h^2\right).
    \end{equation*}
    Note that this estimate is independent of $\mathcal{N}_{\star}$. By substituting this estimate together with \eqref{eq-fm9}--\eqref{eq-fm10} into \eqref{eq-fm6}, we yield
   \begin{equation*}
        \bigl\| \delta_t \xi_u^{n+1} \bigr\| \leq C_4 \tau \sum_{k=2}^{n}\left(\norm{\delta_t\xi_u^{k+1}} + \norm{\xi_u^{k+1}} + \bigl\| \xi_\psi^{k-\frac{3}{2}}\bigr\|  \right)+ \widetilde{C}_3 \left(\tau^2+h^2\right).
    \end{equation*}
    After using Lemmas \ref{lem-estimate1}--\ref{lem-estimate2} and the discrete Gr\"{o}nwall inequality again, we have
    \begin{equation}\label{eq-fm20}
       \max_{1\leq k\leq K}\Bigl( \bigl\| \xi_u^{k} \bigr\| + \bigl\|\xi_\psi^{k-\frac{1}{2}} \bigr\| + \bigl\| \delta_t \xi_u^{k} \bigr\|\Bigr) \leq C^{\ast} \left(\tau^2+h^2\right),
    \end{equation}
    where the constant $C^{\ast}$ is related only to $\mathcal{N}$, $T$ and the indicated derivatives of $g$, but independent of $\tau$, $h$ and $\mathcal{N}_\star$. Repeating \eqref{eq-fm18}--\eqref{eq-fm19} with \eqref{eq-fm20} can yield
    \begin{equation*}
        \bigl\| \wh{\xi}_u^{\,k+\frac{1}{2}} \bigr\|_{\infty} 
        \leq C^{\ast}_{\alpha} (\tau^2+h^2),
    \end{equation*}
    where the constant $C^{\ast}_{\alpha}$ is also related to the fractional power $\alpha$. 
    For sufficiently small $\tau$ and $h$ with $\tau > h$, we further deduce that 
    \begin{equation}\label{eq-fmbegin1}
        \begin{aligned}
            \norm{\xi_u^{n+1}}_{\infty} 
            &= \sum_{k=0}^{n}\left(\norm{\xi_u^{k+1}}_{\infty} - \norm{\xi_u^{k}}_{\infty} \right) \leq \sum_{k=0}^{n}\norm{\xi_u^{k+1}+\xi_u^{k}}_{\infty} = 2\sum_{k=0}^{n} \bigl\|\wh{\xi}_u^{\,k+\frac{1}{2}}\bigr\|_{\infty}\\
            &\leq C^{\ast}_{\alpha} T\tau^{-1} \left(\tau^2+h^2\right) \leq C^{\ast}_{\alpha}T \left(\tau+h\right) \leq 1 - \varepsilon, \quad  0\leq n\leq K-1.
        \end{aligned}
    \end{equation}
    Meanwhile, we can also derive from \eqref{eq-auxerreq3} that 
    \begin{equation*}
        \begin{aligned}
            \bigl\| \xi_\psi^{n+\frac{1}{2}} \bigr\|_{\infty} 
            &\leq \bigl\| \xi_\psi^{\frac{1}{2}} \bigr\|_{\infty} + \sum_{k=1}^{n} \bigl\| \xi_\psi^{k+\frac{1}{2}}+\xi_\psi^{k-\frac{1}{2}}\bigr\|_{\infty} \\
            &\leq \bigl\| \xi_\psi^{\frac{1}{2}} \bigr\|_{\infty} + \sum_{k=1}^{n} \left[2g_\sigma^1 \left(\norm{u^k}_{\infty}+\norm{U_\lambda^k}_{\infty}\right)\norm{\xi_u^{k}}_{\infty} + \norm{\mathtt{r}_\psi^k}_{\infty} \right] \\
            &\leq \widetilde{C}_0(\tau+h^2) + g_\sigma^1\left(4\mathcal{N}+2\right) C^{\ast}_{\alpha}T\sum_{k=1}^{n} \left(\tau+h\right) + \widetilde{C}T\tau,
        \end{aligned}
    \end{equation*}
    where \eqref{eq-mvarbd} and Lemmas \ref{lem-gfunLip}, \ref{lem-truncerror} and \ref{lem-halfrexest} are used. Furthermore, we get
    \begin{equation}
        \bigl\| \xi_\psi^{n+\frac{1}{2}} \bigr\|_{\infty} \leq C^{\dagger}_{\alpha} \tau^{-1} \left(\tau+h\right)  \leq 2C^{\dagger}_{\alpha} = \mathcal{N}_\star - \varepsilon,
    \end{equation}
    where the constant $C^{\dagger}_{\alpha}$  depends only on $C^{\ast}_{\alpha}$, $\mathcal{N}$, $T$, $g_\sigma^1$, $\widetilde{C}_0$ and $\widetilde{C}$, but is independent of $\tau$ and $h$. Here we fix  the constant $\mathcal{N}_\star = 2C^{\dagger}_{\alpha}+\varepsilon$.

    In the case of $\tau \leq h$, for sufficiently small $\tau$ and $h$, it follows from Theorem \ref{thm-1} and the inverse inequality that
    \begin{equation}\label{eq-fmend2}
        \begin{aligned}
            &\norm{\xi_u^k}_{\infty} \leq C_{\mathrm{inv}} h^{-1}\norm{\xi_u^k}\leq C_{\mathrm{inv}}C_3^\star h^{-1}(\tau^2+h^2) \leq C_{\mathrm{inv}}C_3^\star \left(\tau + h\right) \leq 1-\varepsilon,\\
            &\bigl\|\xi_\psi^{k-\frac{1}{2}}\bigr\|_{\infty} \leq C_{\mathrm{inv}}h^{-1}\bigl\|\xi_\psi^{k-\frac{1}{2}}\bigr\| \leq C_{\mathrm{inv}}C_3^\star\left(\tau+h\right) \leq \mathcal{N}_\star-\varepsilon.
        \end{aligned}
    \end{equation}

    Finally, combining  \eqref{eq-fmbegin1}--\eqref{eq-fmend2}, we establish \eqref{eq-uniformbd} without any mesh ratio constraint. 
\end{proof}

In the following, we transfer the modified estimates to the RCND scheme.

\begin{theorem}\label{thm-3}
    Let $e_u^{k}\in \mathsf{C}_h^0$ and  $e_\psi^{k-1/2} \in \mathsf{R}_h^0$ satisfy the original error equations \eqref{eq-orierreq4}--\eqref{eq-orierreq6}. Under the regularity assumptions of Theorem~\ref{thm-1}, there exist a constant $C>0$, independent of $\tau$ and $h$, such that
    \begin{equation}
        \max_{1\leq k\leq K} \bigl\|e_u^{k}\bigr\|_{H^{\frac{\alpha}{2}}_h} \leq C\left(\tau^2+h^2\right), \quad \max_{1\leq k\leq K}\bigl\|e_\psi^{k-\frac{1}{2}}\bigr\|  \leq C \left(\tau^2+h^2\right).
    \end{equation}
\end{theorem}

\begin{proof}
    Based on the uniform bounds established in Theorem \ref{thm-2}, we immediately obtain $\lambda_u = \lambda_\psi \equiv 1$ from \eqref{eq-cutoffcoe}. Therefore, the modified error system \eqref{eq-auxerreq1}--\eqref{eq-auxerreq3} automatically reduces to the original error system \eqref{eq-orierreq4}--\eqref{eq-orierreq6}. Hence $\xi_u^k=e_u^k$ and $\xi_\psi^{k-1/2}=e_\psi^{k-1/2}$, and it can be concluded from Theorem \ref{thm-1} that
    \begin{equation}\label{eq-fm21}
        \max_{1\leq k\leq K}\Bigl( \bigl\| e_u^{k} \bigr\|+\bigl\|e_\psi^{k-\frac{1}{2}}\bigr\| +\bigl\|\delta_t e_u^{k}\bigr\|\Bigr) \leq C \left(\tau^2+h^2\right).
    \end{equation}
    Meanwhile, we also have
    \begin{equation*}
        \max_{1\leq k\leq K}\norm{e_u^k}_{\infty} \leq 1-\varepsilon, \qquad 
        \max_{1\leq k\leq K}\bigl\|e_\psi^{k-\frac{1}{2}}\bigr\|_{\infty} \leq \mathcal{N}_\star-\varepsilon,
    \end{equation*}
    which implies that numerical solutions to the RCND scheme are uniformly bounded, i.e., 
    % it can be seen from \eqref{eq-exsolbound} that 
    \begin{equation}\label{eq-fm22}
        \max_{1\leq k\leq K}\norm{U^k}_{\infty} \leq \mathcal{N} + 1, \qquad 
        \max_{1\leq k\leq K}\bigl\|\Psi^{k-\frac{1}{2}}\bigr\|_{\infty} \leq \mathcal{N} + \mathcal{N}_{\star}.
    \end{equation}

    In addition, by taking an inner product of \eqref{eq-orierreq4} with $\delta_t e_u^{k+1}$, and applying Cauchy-Schwarz inequality and Young's inequality to its real part, we have 
    \begin{equation*}
        \begin{aligned}
            % -\norm{\mathcal{L}^{\alpha} e_u^{k+1}}^2+\norm{\mathcal{L}^{\alpha} e_u^{k}}^2 
            \inner{(-\Delta_h)^{\frac{\alpha}{2}}e_u^{k+1}}{e_u^{k+1}}-\inner{(-\Delta_h)^{\frac{\alpha}{2}}e_u^{k}}{e_u^{k}}
            \leq C\tau \biggl(\sum_{m=1}^{3} \bigl\|\mathrm{F}_{u,m}^{k+\frac{1}{2}}\bigr\|^2 + \norm{\delta_t e_u^{k+1}}^2\biggr).
        \end{aligned}
    \end{equation*}
    Summing the above equation over $k$ from $0$ to $n$ $(0\leq n\leq K-1)$, and using Lemmas \ref{lem:norm_equivalence} and \ref{lem-truncerror}, \eqref{eq-fm21} and \eqref{eq-fm22}, we get
    \begin{equation*}
        \begin{aligned}
            \left(\frac2\pi\right)^\alpha\left|e_u^{n+1}\right|_{H_h^{\frac{\alpha}{2}}}^2
            \leq\inner{(-\Delta_h)^{\alpha/2}e_u^{n+1}}{e_u^{n+1}}
            \leq C\tau\sum_{k=0}^n\left(\sum_{m=1}^3\|\mathrm F_{u,m}^{k+1/2}\|^2+\|\delta_te_u^{k+1}\|^2\right) 
            \leq C\left(\tau^2+h^2\right)^2,
        \end{aligned}
    \end{equation*}
    which, together with \eqref{eq-fm21}, yields that 
    \begin{equation*}
        \norm{e_u^{n+1}}_{H^{\frac{\alpha}{2}}_h} = \sqrt{\norm{e_u^{n+1}}^2 + \left|e_u^{n+1}\right|_{H^{\frac{\alpha}{2}}_h}^2}\leq C\left(\tau^2+h^2\right).
    \end{equation*}
    This completes the proof of the unconditional optimal-order convergence for the RCND scheme.
\end{proof}

%\begin{remark}
%    Consider the equation on $\Omega=(a,b)^3$ and set the discrete inner products and norms in the 3D space.
%    Note that  all lemmas Section \ref{sec-preliminaries}, except the Lemma \ref{lem-uniform_maxest}, are applicable to 3D equations. If $3/2<\alpha\leq2$,  then Lemma \ref{lem-uniform_maxest} can be extended to 3D cases. Specifically, the bound in inequality \eqref{eq-lowfrees1} becomes $CR^3\|w\|_1^2$, which then yields $\|w\|^{2+2\alpha/3}\leq C E(w)\|w\|_1^{2\alpha/3}$, corresponding to \eqref{eq-estf2}. Subsequently, a similar analysis process is carried out, the second bound of \eqref{eq-resolvemaxes1} requires $\int_0^\infty e^{-t}t^{-3/(2\alpha)}\,\rmd t<\infty$, namely $\alpha > 3/2$. Thus, by employing the same arguments as in Section \ref{sec-scheme}--\ref{sec-convergence}, we obtain the unconditional optimal-order convergence for 3D spaces, where the 3D inverse inequality $\|v\|_\infty\leq C_{\mathrm{inv}}h^{-3/2}\|v\|$ is used.
%\end{remark}
%

\section{Extension to three dimensions}\label{sec-three-dimensions}

Consider the SFNS equation on $\Omega=(a,b)^3$ with homogeneous exterior data. Use a uniform cubic grid with spacing $h$, and extend the grid function by zero outside $\Omega$. Then, the discrete inner product is defined as $\langle v,w\rangle=h^3\sum_{\mathbf{j}\in\mathbb{Z}^3}v_{\mathbf{j}}\overline{w}_{\mathbf{j}}$, with norms $\|v\|^2=\sqrt{\langle v,v\rangle}$ and $\|v\|_1=h^3\sum_{\mathbf{j}\in\mathbb{Z}^3}|v_{\mathbf j}|$. Meanwhile, the Fourier definition of $H_h^s$ employs the interval $[-\pi/h,\pi/h]^3$ and the factor $(2\pi)^{-3}$. According to \cite{hao2021fractional}, replace the spatial operator by
\begin{equation}\label{eq-3d-operator}
    A_hv_{\mathbf{j}}=\frac{1}{h^{\alpha}}\sum_{\mathbf{m}\in\mathbb Z^3}a_{\mathbf{m}}^{(\alpha,3)}v_{\mathbf{j}+\mathbf{m}}, 
\end{equation}
where the coefficients are given by
\begin{equation}
    a_{\mathbf{m}}^{(\alpha,3)}=\frac{1}{(2\pi)^3}\int_{[-\pi,\pi]^3}\left(4\sum_{\ell=1}^3\sin^2\frac{\theta_\ell}{2}\right)^{\frac{\alpha}{2}}e^{-\bi\mathbf{m}\cdot\theta} \,\rmd\theta.
\end{equation}
The symbol retains the coupling between coordinate directions for $1<\alpha<2$; when $\alpha=2$, $A_h$ is the standard seven-point negative Laplacian. The half-step initialization and all time updates in \eqref{eq-numsma}--\eqref{eq-numsmd} keep their form, with $(-\Delta_h)^{\alpha/2}$ replaced by $A_h$.

The coefficient signs, symmetry, and fractional semi-norm equivalence can be extended by the same arguments. Hence, the unique solvability in Theorem \ref{thm-uniquesolve} and the discrete conservation laws in Theorem \ref{thm-conservation} remain valid. The consistency estimates in Lemma \ref{lem-truncerror} also extend, since its proof relies solely on the same operator properties and pointwise nonlinear estimates.

The main changes related to dimensions appear in Lemma \ref{lem-uniform_maxest}. Specifically, the bound in inequality \eqref{eq-lowfrees1} becomes $CR^3\|w\|_1^2$, which then yields the inequality $\|w\|^{2+2\alpha/3}\leq C E(w)\|w\|_1^{2\alpha/3}$, corresponding to \eqref{eq-estf2}. Using the same notation $P(t)=e^{-th^{-\alpha}A_h}$, and analogous to \eqref{eq-maxnorm-es1}, we obtain 
\begin{equation*}
    \norm{P(t)f} \leq C t^{-3/(2\alpha)} \norm{f},
\end{equation*}
which leads to the second bound of \eqref{eq-resolvemaxes1} requires $\int_0^\infty e^{-t}t^{-3/(2\alpha)}\,\rmd t<\infty$, namely $\alpha>3/2$. Therefore, the conclusion of Lemma \ref{lem-uniform_maxest} in the three-dimensional case is as follows:
\begin{equation}\label{eq-3d-maxnorm}
 \|v\|_\infty\leq C_\alpha\bigl(\|v\|+\|A_hv\|\bigr), \qquad 3/2< \alpha \leq 2.
\end{equation}

Subsequently, by employing the same arguments as in Section \ref{sec-convergence}, we can extend the result of unconditional optimal-order convergence to the three-dimensional case. The main difference begins with Theorem \ref{thm-2}.
For $\tau>h$, \eqref{eq-3d-maxnorm} gives the same midpoint bounds; the time-difference estimates and the choice of the cut-off level remain applicable as well. For $\tau\leq h$, the inverse inequality becomes $\|v\|_\infty \leq C_{\mathrm{inv}}h^{-3/2}\|v\|$. Thus, for either $\xi=\xi_u^k$ or $\xi=\xi_\psi^{k-1/2}$,
\begin{equation*}
    \|\xi\|_\infty\leq C_{\mathrm{inv}}C_{\star} h^{-3/2}(\tau^2+h^2) \leq 2C_{\mathrm{inv}}C_\star h^{1/2}\longrightarrow 0.
\end{equation*}
After fixing the cut-off level as in the two-dimensional proof, sufficiently small $\tau$ and $h$, independently of each other, ensure that both cut-off coefficients equal one. For $3/2<\alpha\leq2$, based on the three-dimensional fractional semi-norm equivalence, the energy argument of Theorem \ref{thm-3} then yields
\begin{equation}\label{eq-3d-convergence}
 \max_{1\leq k\leq K}\|u^k-U^k\|_{H_h^{\alpha/2}}
 +\max_{1\leq k\leq K}\|\psi^{k-1/2}-\Psi^{k-1/2}\|
 \leq C(\tau^2+h^2),
\end{equation}
for sufficiently small $\tau$ and $h$, with $C$ independent of both. 

In a word, the three-dimensional scheme has the same optimal second-order convergence without a mesh-ratio restriction, but the range of $\alpha$ is limited to $3/2<\alpha\leq 2$.

\section{Numerical experiments}\label{sec-numerics}

In this section, we provide several numerical experiments to verify the theoretical analysis of the proposed RCND scheme \eqref{eq-numsma}--\eqref{eq-numsmd}. 
For the discrete fractional Laplacian $(-\Delta_h)^{\frac{\alpha}{2}}$, the matrix coefficients $a_{i,j}^{(\alpha)}$ given in \eqref{eq-compute-intcoeff} can be approximated by the trapezoidal rule, i.e., 
\begin{equation}\label{eq-coeffcompute}
    a_{i,j}^{(\alpha)} \approx \frac{1}{L^2} \sum_{\zeta_1=0}^{L-1}\sum_{\zeta_2=0}^{L-1}\left[4 \sin^2\left(\frac{\zeta_1\delta}{2}\right) + 4 \sin^2\left(\frac{\zeta_2\delta}{2}\right)\right]^{\frac{\alpha}{2}} \! e^{-\bi(i\zeta_1 + j\zeta_2)\delta}, \quad 0\leq i,j\leq M-2, 
\end{equation}
where $\delta = 2\pi/L$ and $L\gg M$. Obviously, we can use FFT algorithm to evaluate the coefficient $a_{i,j}^{(\alpha)}$ efficiently through the above formula. 
In the following examples, we take $L=2^{15}$.
Thus, we employ the Bi-CGSTAB iterative method \cite{van1992bi} to solve the resulting algebraic system, with the iteration tolerance set to $\mathrm{Tol}=10^{-16}$.

Next, we introduce the  error measures. Set $e^K(\tau,h) = U^K(\tau,h) - u^K(\tau,h)$, where $u^K$ and $U^K$ denote the exact solution and numerical solution at $t=T$, respectively, and $(\tau,h)$ denotes the corresponding space-time step.
For the primary variable, we define the errors  
\begin{equation*}
  \mathrm{Err}_{u}^{2} = \bigl\|e^K(\tau,h) \bigr\|, \qquad \mathrm{Err}_{u}^{\alpha} = \bigl\|e^K(\tau,h) \bigr\|_{\mathcal{A}^{\alpha}}.
\end{equation*}
 Similarly, the error for the relaxation variable at $t=T-\tau/2$ is denoted by 
\begin{equation*}
    \mathrm{Err}_{\psi}^{2} = \bigl\|\Psi^{K-\frac{1}{2}}(\tau,h) - \psi^{K-\frac{1}{2}}(\tau,h) \bigr\|.
\end{equation*}
The superscript $\alpha$ in $\mathrm{Err}_{u}^{\alpha}$ denotes the following discrete energy norm,  
% It is an equivalent norm, rather than a quadrature approximation to the Fourier-defined Sobolev norm:
\begin{equation}\label{eq-appro-sobnorm}
    \|e^K(\tau,h)\|_{\mathcal A^\alpha}= \sqrt{\bigl\|e^K(\tau,h)\bigr\|^2 + \inner{(-\Delta_h)^{\frac{\alpha}{2}}e^K(\tau,h)}{e^K(\tau,h)}},
\end{equation}
which is uniformly equivalent to the $H_h^{\alpha/2}$-norm according to the Lemma \ref{lem:norm_equivalence}. Consequently, we use \eqref{eq-appro-sobnorm} to measure $\|e^K(\tau,h)\|_{H^{\alpha/2}_h}$.

\begin{exm}[Manufactured Solution]\label{exm-1}
    Consider the following SFNS problem 
    \begin{equation}
      \begin{aligned}
        \bi u_t - (-\Delta)^{\frac{\alpha}{2}}u + g \left(|u|^2\right)u &= f,  \quad (x,y) \in \Omega,&& \quad 0<t<T, \\
        u &= 0, \quad (x,y)\in \Omega^c,&& \quad 0<t<T,
      \end{aligned}
    \end{equation}
    where $\Omega=(-4,4)^2$, $T=1$. Moreover, we set the exact solution as $u(x,y,t) = \exp(-\bi t)\exp(-x^2-y^2)$, and choose $g(s) = s, s^2, s^3$.
\end{exm}

Since it is difficult to derive the explicit expression of the source function $f(x,y,t)$, we compute it approximately with a finer space grid size, i.e., $f\approx \bi u_t - (-\Delta_{h})^{\frac{\alpha}{2}}u + g \left(|u|^2\right)u$ with $h=1/480$.
Tables \ref{tab:exm1-tsorder-a}, \ref{tab:exm1-tsorder-b} and \ref{tab:exm1-tsorder-c} correspond to the functions $g(s) = s, s^2, s^3$, respectively. 
In these tables, for the different $\alpha$, we list the $L^2$-norm and energy-norm errors for the primary variable $U$, and the $L^2$-norm  errors for the relaxation variable $\Psi$. 
The reported rates are approximately second order over the tested meshes, which is consistent with the theoretical results established in Theorem \ref{thm-3}.

\begin{table}[!tbh]%[htbp]
\centering
\caption{Errors and convergence rates for Example \ref{exm-1} with $g(s)=s$.}
\label{tab:exm1-tsorder-a}
\begin{tabular}{cccccccccc}
\toprule

& $\tau=h$ & $\mathrm{Err}_{u}^{2}$ & $\mathrm{Rate}$ & $\mathrm{Err}_{u}^{\alpha}$ & $\mathrm{Rate}$ & $\mathrm{Err}_{\psi}^{2}$ & $\mathrm{Rate}$ \\
\midrule
\multirow{4}{*}{$\alpha=1.2$} 
& 1/5  & 1.9693e-02 & $-$   & 3.8277e-02 & $-$  & 1.3318e-02 & $-$  \\
 
& 1/10 & 5.0453e-03 & 1.96  & 9.8210e-03 & 1.96 & 3.7285e-03 & 1.84 \\
 
& 1/20 & 1.2676e-03 & 1.99  & 2.4709e-03 & 1.99 & 9.8000e-04 & 1.93  \\

& 1/40 & 3.1787e-04 & 2.00  & 6.2068e-04 & 1.99 & 2.5371e-04 & 1.95  \\
\midrule
\multirow{4}{*}{$\alpha=1.5$}
& 1/5  & 2.6464e-02 & $-$   & 5.4390e-02 & $-$  & 2.2794e-02 & $-$ \\
 
& 1/10 & 6.6368e-03 & 2.00  & 1.3503e-02 & 2.01 & 5.4626e-03 & 2.06  \\
 
& 1/20 & 1.6580e-03 & 2.00  & 3.3662e-03 & 2.00 & 1.3806e-03 & 1.98 \\

& 1/40 & 4.1239e-04 & 2.00  & 8.3783e-04 & 2.00 & 3.4703e-04 & 1.99 \\
\midrule
\multirow{4}{*}{$\alpha=1.8$}
& 1/5  & 3.0207e-02 & $-$   & 6.3375e-02 & $-$  & 2.6044e-02 &  $-$\\
 
& 1/10 & 7.4382e-03 & 2.02  & 1.5557e-02 & 2.03 & 5.9096e-03 & 2.14  \\
 
& 1/20 & 1.8527e-03 & 2.01  & 3.8732e-03 & 2.01 & 1.4970e-03 & 1.98 \\

& 1/40 & 4.6062e-04 & 2.01  & 9.6295e-04 & 2.01 & 3.7478e-04 & 2.00 \\
\bottomrule 
\end{tabular}
\end{table}

\begin{table}[!tbh]
\centering
\caption{Errors and convergence rates for Example \ref{exm-1} with $g(s)=s^2$.}
\label{tab:exm1-tsorder-b}
\begin{tabular}{cccccccccc}
\toprule

& $\tau=h$ & $\mathrm{Err}_{u}^{2}$ & $\mathrm{Rate}$ & $\mathrm{Err}_{u}^{\alpha}$ & $\mathrm{Rate}$ & $\mathrm{Err}_{\psi}^{2}$ & $\mathrm{Rate}$ \\
\midrule
\multirow{4}{*}{$\alpha=1.2$} 
& 1/5  & 2.3551e-02 & $-$   & 4.6463e-02 & $-$  & 2.9326e-02 & $-$  \\
 
& 1/10 & 6.1436e-03 & 1.94  & 1.2145e-02 & 1.94 & 8.5879e-03 & 1.77 \\
 
& 1/20 & 1.5509e-03 & 1.99  & 3.0693e-03 & 1.98 & 2.2787e-03 & 1.91  \\

& 1/40 & 3.8872e-04 & 2.00  & 7.7178e-04 & 1.99 & 5.9101e-04 & 1.95  \\
\midrule
\multirow{4}{*}{$\alpha=1.5$}
& 1/5  & 3.2627e-02 & $-$   & 6.7321e-02 & $-$  & 4.9492e-02 & $-$ \\
 
& 1/10 & 8.2091e-03 & 1.99  & 1.6725e-02 & 2.01 & 1.1996e-02 & 2.04  \\
 
& 1/20 & 2.0521e-03 & 2.00  & 4.1706e-03 & 2.00 & 3.0397e-03 & 1.98 \\

& 1/40 & 5.1050e-04 & 2.01  & 1.0381e-03 & 2.01 & 7.6520e-04 & 1.99 \\
\midrule
\multirow{4}{*}{$\alpha=1.8$}
& 1/5  & 3.6992e-02 & $-$   & 7.7007e-02 & $-$  & 5.3764e-02 &  $-$\\
 
& 1/10 & 9.0293e-03 & 2.03  & 1.8707e-02 & 2.04 & 1.2218e-02 & 2.14  \\
 
& 1/20 & 2.2452e-03 & 2.01  & 4.6466e-03 & 2.01 & 3.0804e-03 & 1.99 \\

& 1/40 & 5.5798e-04 & 2.01  & 1.1546e-03 & 2.01 & 7.7107e-04 & 2.00 \\
\bottomrule 
\end{tabular}
\end{table}

\begin{table}[!tbh]
\centering
\caption{Errors and convergence rates for Example \ref{exm-1} with $g(s)=s^3$.}
\label{tab:exm1-tsorder-c}
\begin{tabular}{cccccccccc}
\toprule

& $\tau=h$ & $\mathrm{Err}_{u}^{2}$ & $\mathrm{Rate}$ & $\mathrm{Err}_{u}^{\alpha}$ & $\mathrm{Rate}$ & $\mathrm{Err}_{\psi}^{2}$ & $\mathrm{Rate}$ \\
\midrule
\multirow{4}{*}{$\alpha=1.2$} 
& 1/5  & 2.8048e-02 & $-$   & 5.7157e-02 & $-$  & 5.1980e-02 & $-$  \\
 
& 1/10 & 7.4850e-03 & 1.91  & 1.5304e-02 & 1.90 & 1.5960e-02 & 1.70 \\
 
& 1/20 & 1.8997e-03 & 1.98  & 3.8857e-03 & 1.98 & 4.2706e-03 & 1.90  \\

& 1/40 & 4.7662e-04 & 1.99  & 9.7819e-04 & 1.99 & 1.1087e-03 & 1.95  \\
\midrule
\multirow{4}{*}{$\alpha=1.5$}
& 1/5  & 3.9481e-02 & $-$   & 8.3357e-02 & $-$  & 8.6631e-02 & $-$ \\
 
& 1/10 & 9.9024e-03 & 2.00  & 2.0522e-02 & 2.02 & 2.0894e-02 & 2.05  \\
 
& 1/20 & 2.4710e-03 & 2.00  & 5.1012e-03 & 2.01 & 5.2663e-03 & 1.99 \\

& 1/40 & 6.1448e-04 & 2.01  & 1.2687e-03 & 2.01 & 1.3249e-03 & 1.99 \\
\midrule
\multirow{4}{*}{$\alpha=1.8$}
& 1/5  & 4.3693e-02 & $-$   & 9.1624e-02 & $-$  & 8.8743e-02 &  $-$\\
 
& 1/10 & 1.0468e-02 & 2.06  & 2.1760e-02 & 2.07 & 1.9546e-02 & 2.18  \\
 
& 1/20 & 2.5916e-03 & 2.01  & 5.3767e-03 & 2.02 & 4.8693e-03 & 2.01 \\

& 1/40 & 6.4338e-04 & 2.01  & 1.3343e-03 & 2.01 & 1.2161e-03 & 2.00 \\
\bottomrule 
\end{tabular}
\end{table}

\begin{exm}[Conservation Verification]\label{exm-2}
    Consider the SFNS equation \eqref{eq-SFNS} with $g(s) = s,s^2,s^3$, and the initial condition is chosen as
    \begin{equation*}
        u(x,y,0) = \exp(-x^2 - y^2), \qquad (x,y) \in \Omega=(-10,10)^2,
    \end{equation*}
    and we impose the homogeneous boundary condition \eqref{eq-bondcd}.
\end{exm}

\begin{figure}[!htbp]
    \centering
    \includegraphics[width=0.27\textwidth]{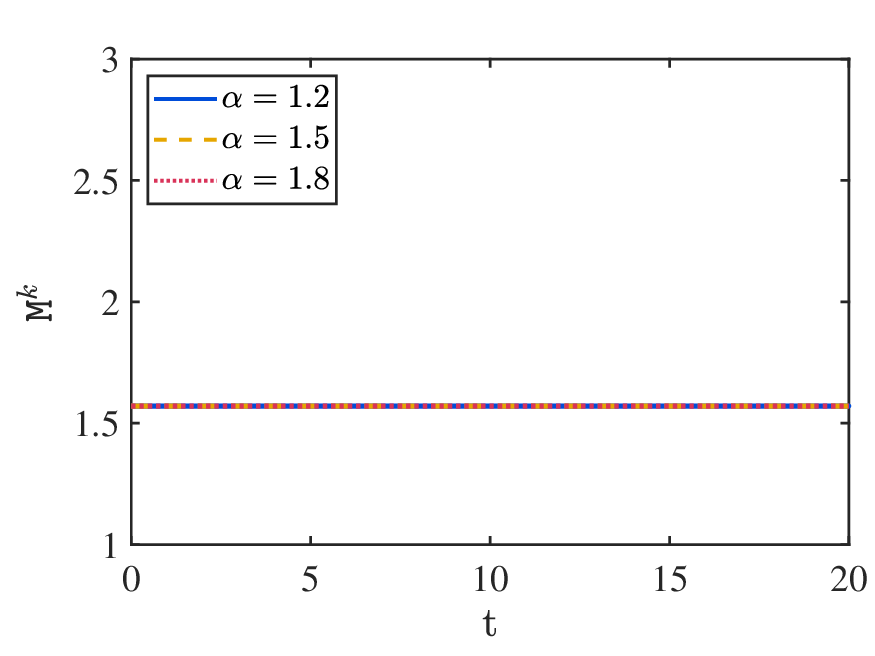}
    \hspace{2mm}
    \includegraphics[width=0.27\textwidth]{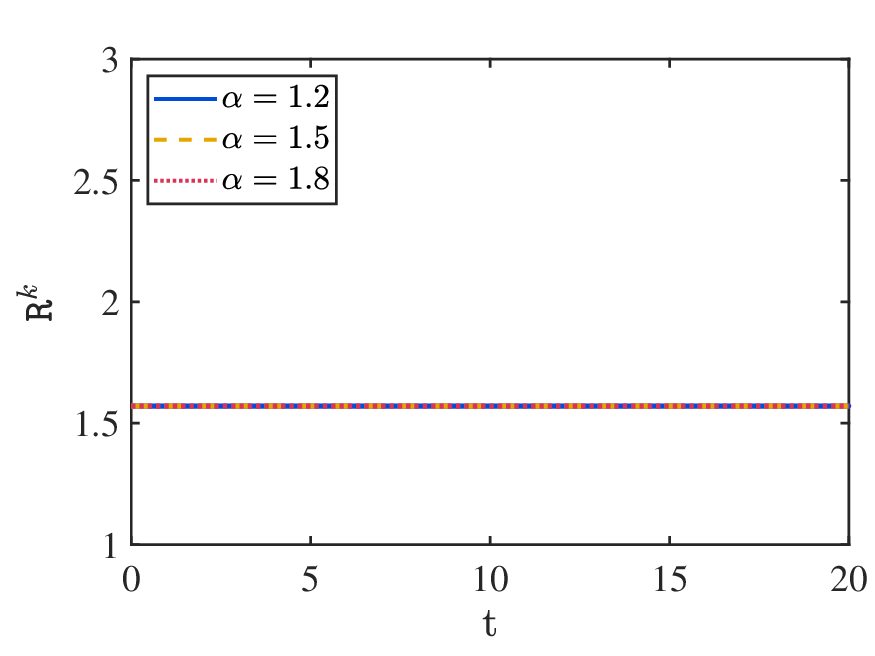}
    \hspace{2mm}
    \includegraphics[width=0.27\textwidth]{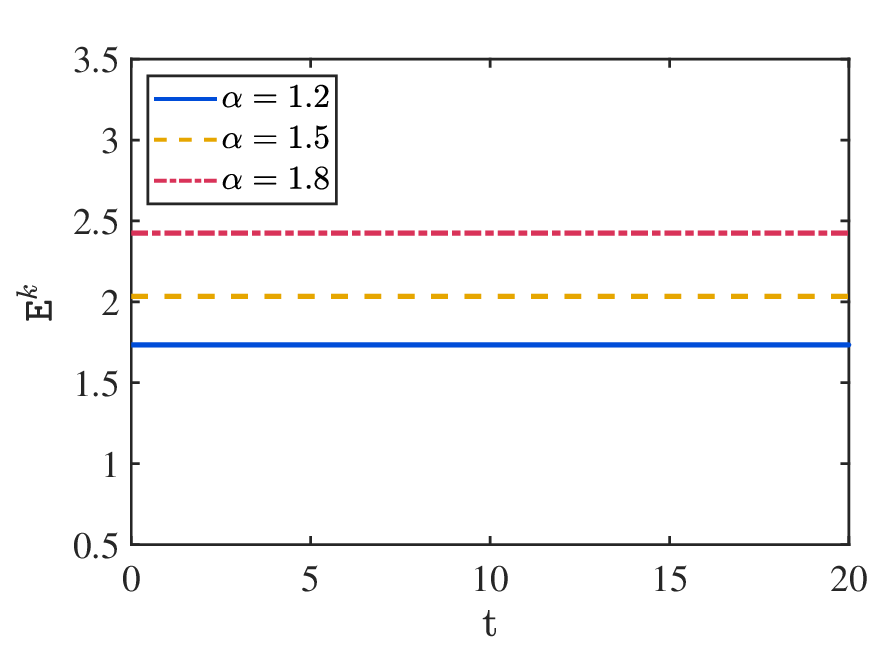}
    
    \par\vspace{2pt}
    \includegraphics[width=0.27\textwidth]{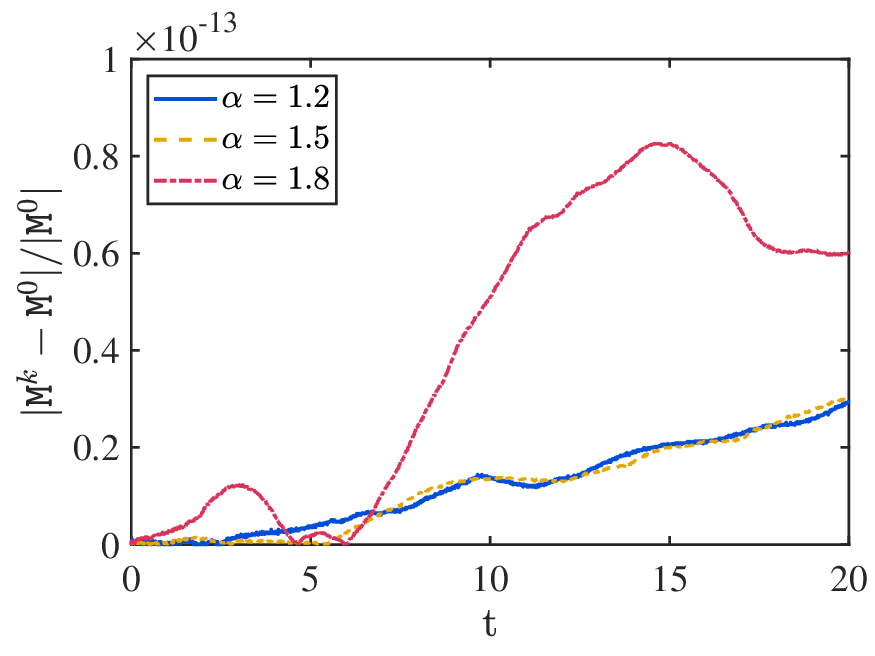}
    \hspace{2mm}
    \includegraphics[width=0.27\textwidth]{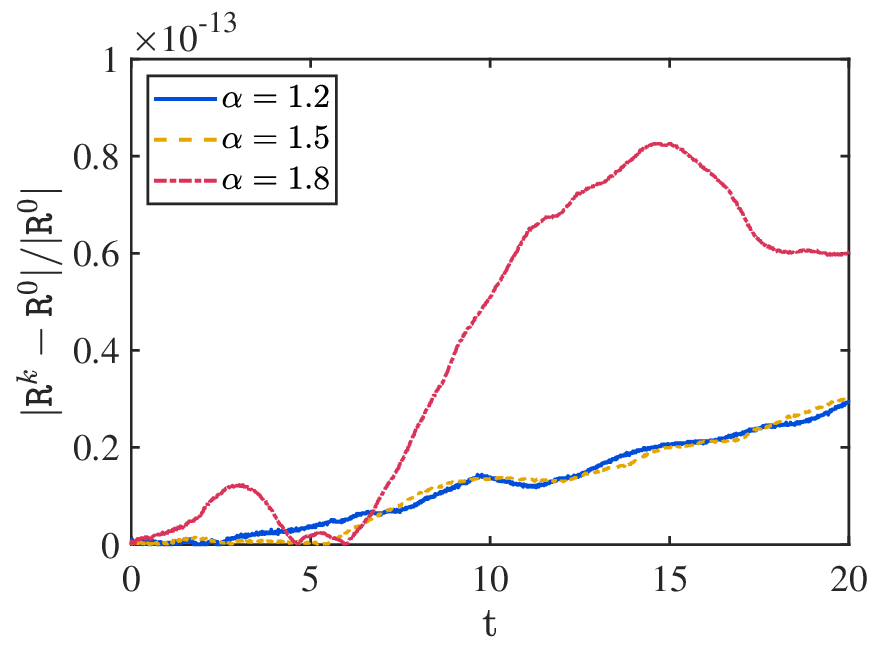}
    \hspace{2mm}
    \includegraphics[width=0.27\textwidth]{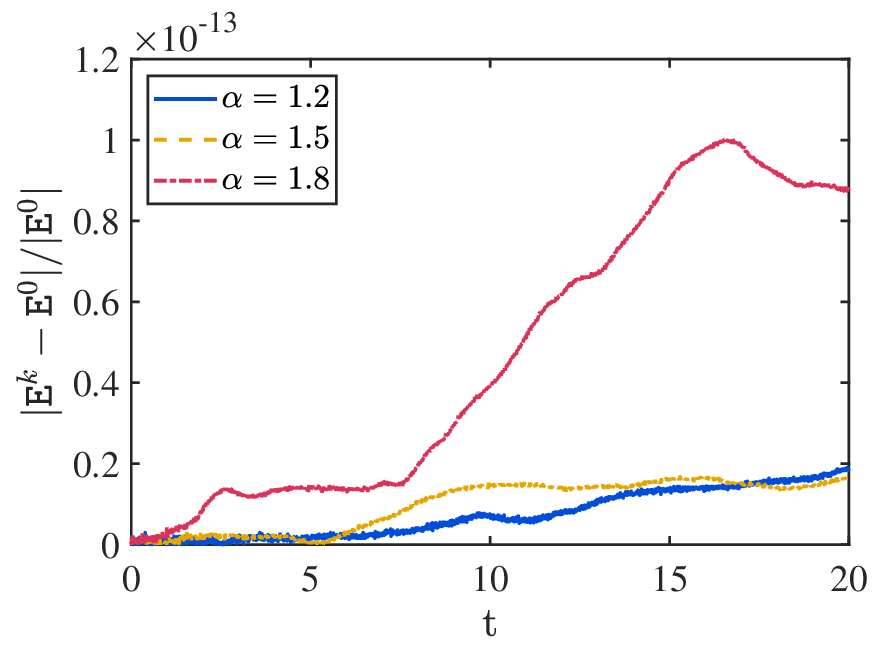}
    \caption{\centering Time evolution of the invariants and their relative errors in Example \ref{exm-2}  for different values of $\alpha$, with $\tau=1/100$, $h=1/20$ and $g(s)=s$.  Left: $\mathtt M^k$; middle: $\mathtt R^k=\mathtt M^k$; right: the modified energy $\mathtt E^k$.}
    \label{fig-exm3a}
\end{figure}

In this example, we verify the conservation laws of Theorem \ref{thm-conservation}. In the cubic case, we derive the discrete mass and modified energy $\mathtt{E}^k$ conservation, where the mass has two forms: $\mathtt{M}^k$ and $\mathtt{R}^k$. Figure \ref{fig-exm3a} plots these conserved quantities and their relative errors, with $\tau=1/100$ and $h=1/20$. For the different values of $\alpha$, the mass $\mathtt{M}^k$, $\mathtt{R}^k$ and the energy $\mathtt{E}^k$ are well preserved over a long time interval, while their relative errors remain on the order of $10^{-13}$. Therefore, the conservation property for the case $g(s)=s$ is confirmed. 

\begin{figure}[!htbp]
    \centering 
    \includegraphics[width=0.28\textwidth]{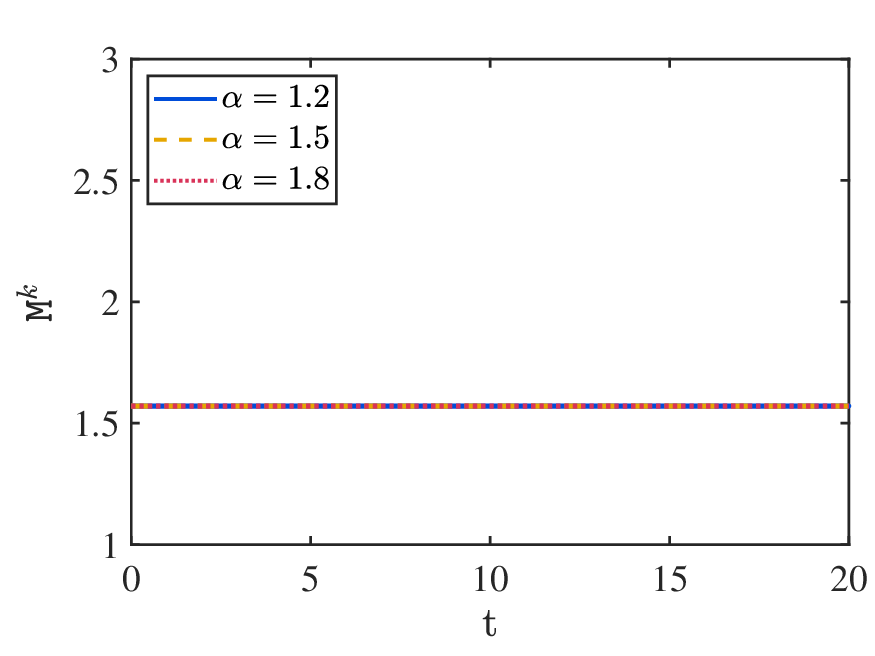}
    \hspace{0.2cm}
    \includegraphics[width=0.28\textwidth]{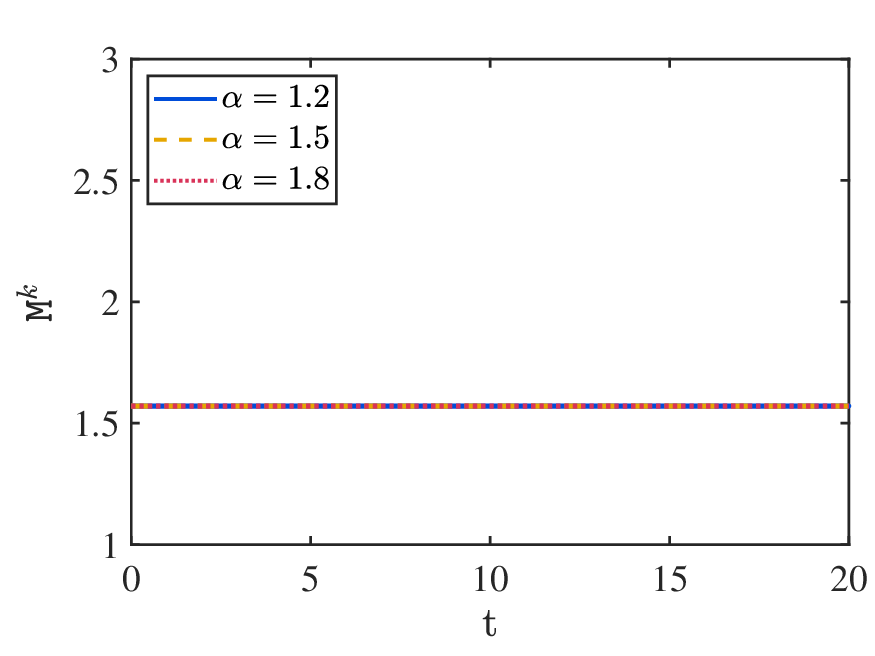}
    
    \par\vspace{2pt}
    \includegraphics[width=0.28\textwidth]{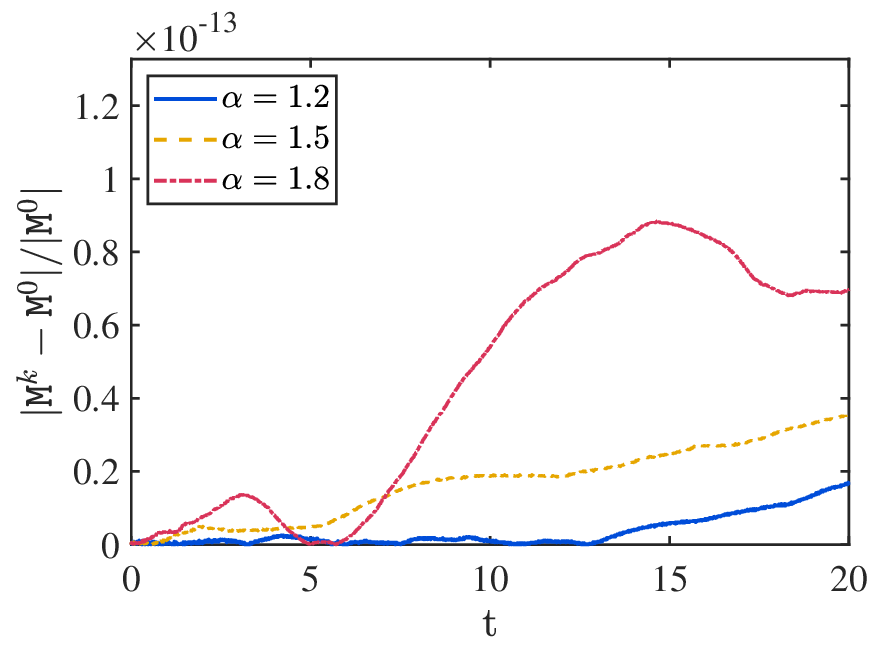}
    \hspace{0.2cm}
    \includegraphics[width=0.28\textwidth]{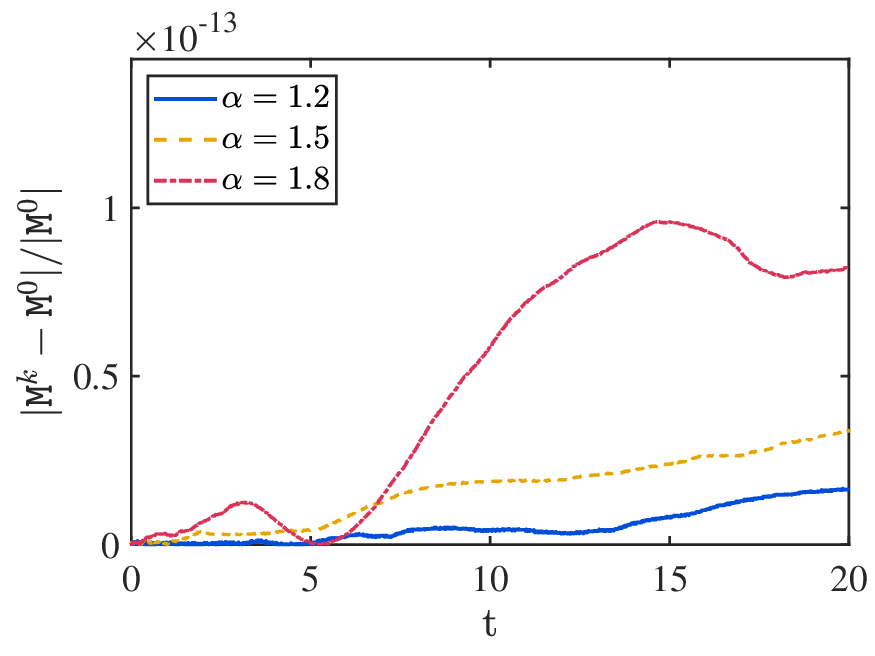}
    \caption{\centering Time evolution of the mass and its relative error in Example \ref{exm-2} for different values of $\alpha$, with $\tau=1/100$ and $h=1/20$.  Left: $g(s)=s^2$; right: $g(s)=s^3$.}
    \label{fig-exm3b}
\end{figure}

For the general nonlinearity, we only obtain the discrete mass conservation from Theorem \ref{thm-conservation}. Here, we set two cases: $g(s) = s^2$ and $g(s) = s^3$.
Under the same spatio-temporal step sizes as Figure \ref{fig-exm3a}, the evolution of the mass $\mathtt{M}^k$ and its relative error over a long time interval is shown in Figure \ref{fig-exm3b}. Throughout the simulation, the relative error for various $\alpha$ remains at the level of $10^{-13}$, which confirms the discrete mass conservation for the general nonlinearity.
Since the initial condition is independent of $\alpha$, the mass should be independent of $\alpha$ as well. Figures \ref{fig-exm3a} and \ref{fig-exm3b} demonstrate this.

\begin{figure}[!htbp]
    \centering
    \includegraphics[width=0.27\textwidth]{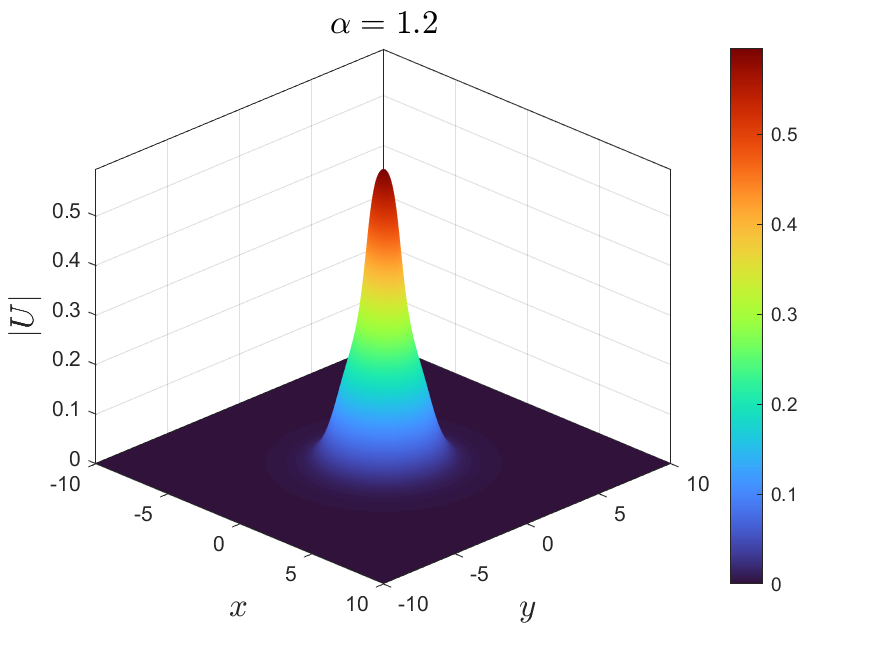}
    \hspace{2mm}
    \includegraphics[width=0.27\textwidth]{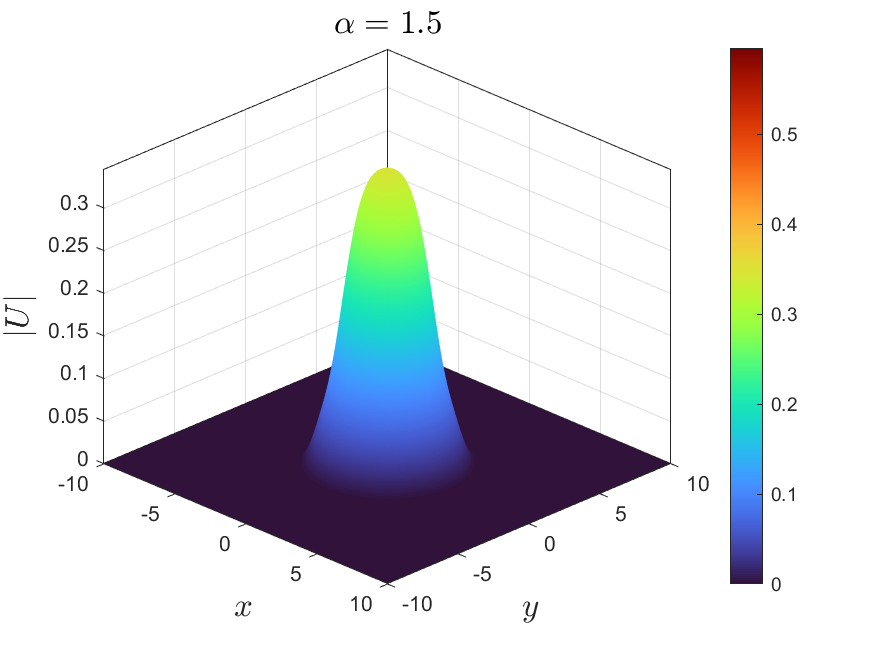}
    \hspace{2mm}
    \includegraphics[width=0.27\textwidth]{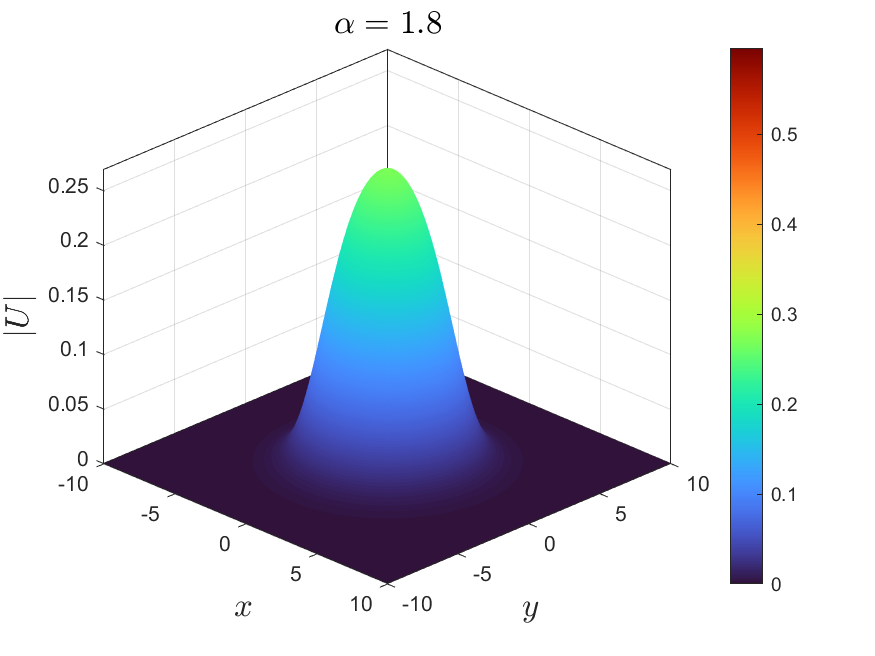}
    
    \par\vspace{2pt}
    \includegraphics[width=0.27\textwidth]{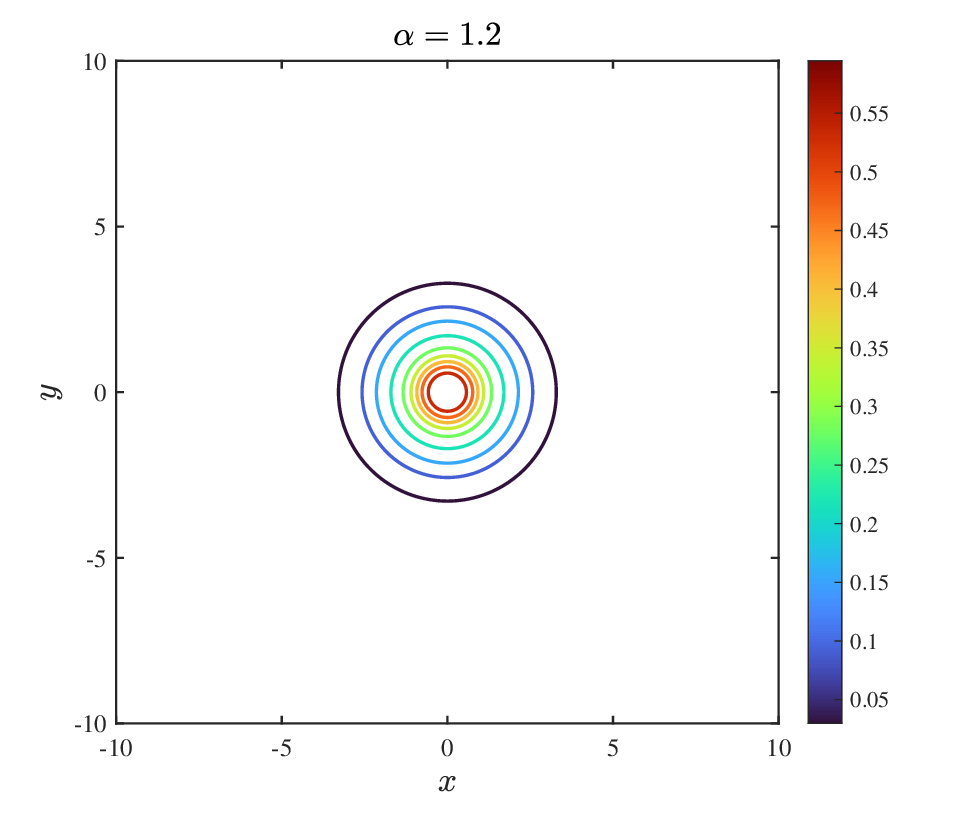}
    \hspace{2mm}
    \includegraphics[width=0.27\textwidth]{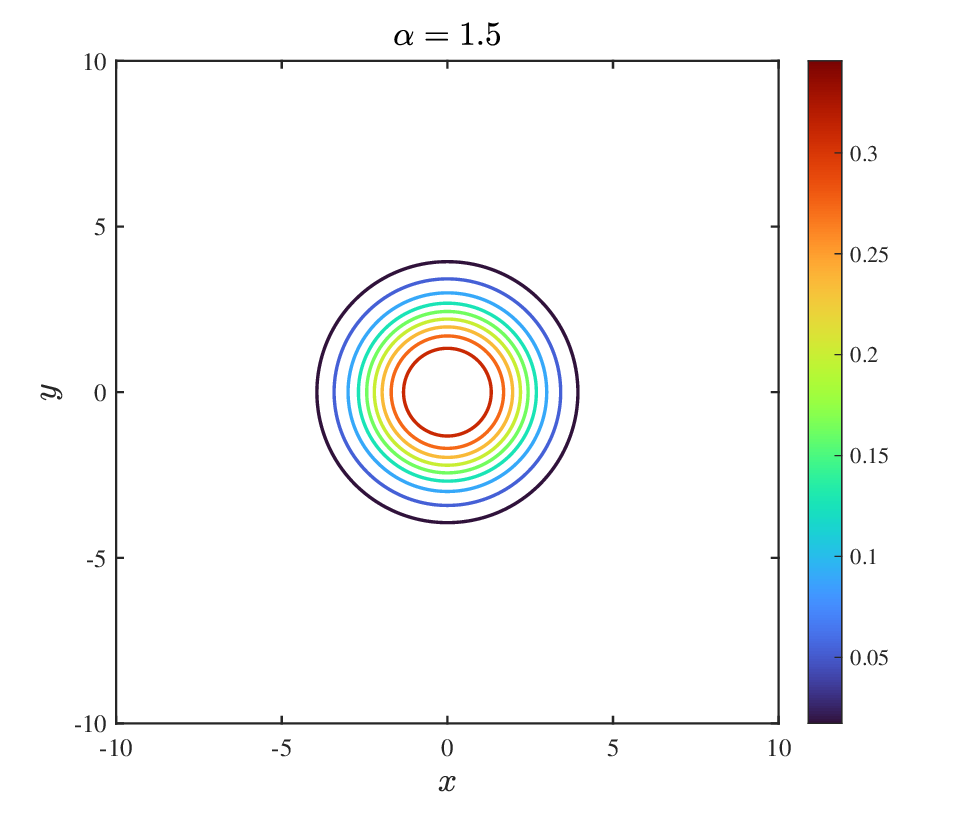}
    \hspace{2mm}
    \includegraphics[width=0.27\textwidth]{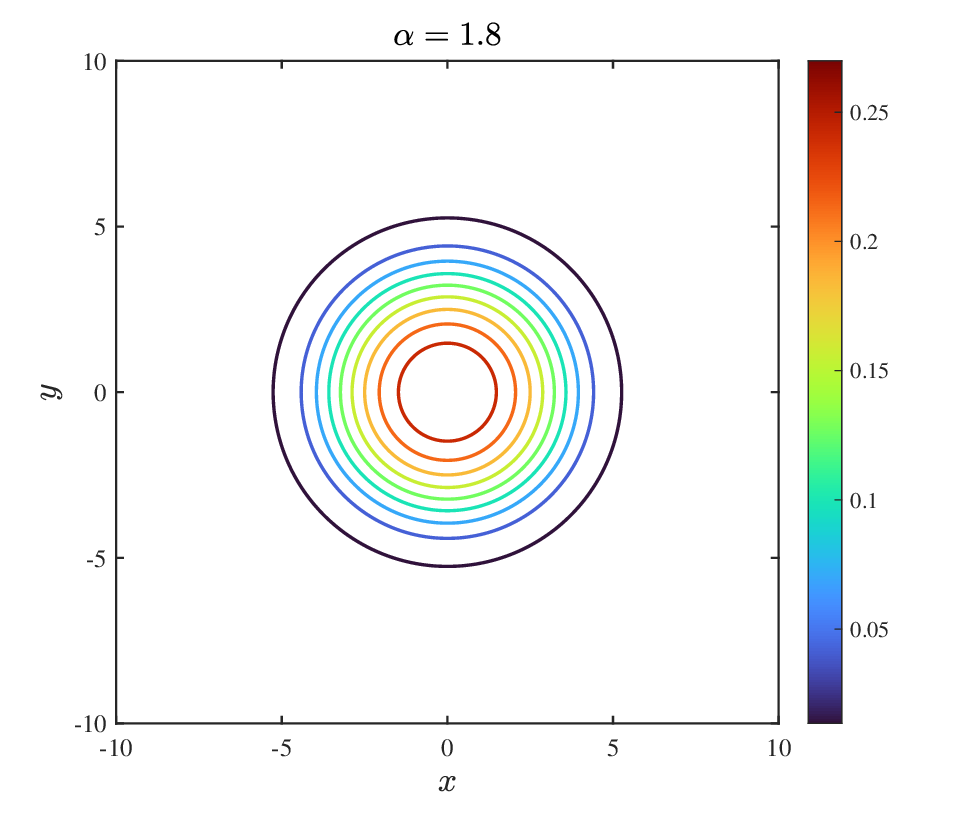}
    \caption{\centering Profile (top) and contour plot (bottom) of the numerical amplitude $|U|$ for various $\alpha$ in Example \ref{exm-2},  with $\tau=1/100$, $h= 1/50$ and $g(s)=s$.}
    \label{fig-exm3c}
\end{figure}

For the cubic case, in which both mass and modified energy are conserved, Figure \ref{fig-exm3c} shows the profiles of the numerical solution $|U|$ and their corresponding contours at $T=1$ for $\alpha=1.2,1.5,1.8$, with $\tau=1/100$ and $h=1/50$. It can be observed from Figure \ref{fig-exm3c} that  the circular wave profile propagates toward the boundary more rapidly as $\alpha$ increases.

% Keep all figures before the conclusion, without forcing an empty figure page.
\FloatBarrier
\section{Conclusion}\label{sec-conclusion}

In this paper, we proposed a relaxation Crank-Nicolson difference scheme for the two-dimensional SFNS equation, where the integral fractional Laplacian is approximated by a fractional central difference scheme based on the generating function. 
To overcome the difficulty caused by the failure of the fractional Sobolev inequality in two-dimensional space, we developed a uniform discrete $L^{\infty}$-norm estimate in Lemma \ref{lem-uniform_maxest}. 
Moreover, we demonstrated that the scheme is uniquely solvable and satisfies the discrete conservation laws. Specifically, it conserves discrete mass for general real-valued $g$, and a modified discrete energy for $g(s)=\beta s$.
By introducing a modified error system using a cut-off tenique, we derived the unconditional optimal-order error estimates of  $\mathcal O(\tau^2+h^2)$ for the primary variables in the $H_h^{\alpha/2}$\!-norm and for the relaxation variables in the $L^2$-norm. Finally, the numerical experiments illustrate the observed second-order rates and the discrete conservation properties.

\section*{Funding}
The research of Dongling Wang is supported in part by the National Natural Science Foundation of China under Grant 12671482, and the Natural Science Foundation of Hunan Province under Grant 2026JJ50363. The work of Lisen Ding is supported by the Postgraduate Scientific Research Innovation Project of Xiangtan University, China (No. XDCX2025Y213).

\section*{Acknowledgements}
During the preparation of this work, we used ChatGPT (OpenAI) to assist with language polishing. The authors reviewed all suggested edits and take full responsibility for the published content.

 \bibliographystyle{unsrt}
\bibliography{references}

\end{document}